\documentclass{article}   	
\usepackage{geometry}               
\usepackage{graphicx}
\usepackage[usenames,dvipsnames]{xcolor}
\usepackage{amssymb,amsmath,amsthm,amsfonts}

\allowdisplaybreaks

\usepackage[english]{babel}
\usepackage[utf8]{inputenc}
\usepackage[colorlinks]{hyperref}
\hypersetup{linkcolor=blue,citecolor=blue,filecolor=black,urlcolor=blue}

\usepackage{mathtools}
\usepackage{dsfont}

\usepackage[sc]{mathpazo}

\usepackage{tikz}
\newtheorem{proposition}{Proposition}[section]
\newtheorem{theorem}[proposition]{Theorem}
\newtheorem{corollary}[proposition]{Corollary}
\newtheorem{lemma}[proposition]{Lemma}

\theoremstyle{definition}
\newtheorem{definition}[proposition]{Definition}
\theoremstyle{remark}
\newtheorem{remark}[proposition]{Remark}
\numberwithin{equation}{section}

\newcommand{\eps}{\varepsilon}

\newcommand{\N}{{\mathbb{N}}}

\newcommand{\R}{{\mathbb{R}}}

\newcommand{\sphere}{{\mathbb{S}}}

\newcommand{\loc}{{\text{loc}}}

\newcommand{\bm}{\mathbf}
\newcommand{\bx}{{\bm{x}}}
\newcommand{\by}{{\bm{y}}}
\newcommand{\bz}{{\bm{z}}}

\newcommand{\bp}{{\bm{p}}}

\newcommand{\Lcal}{{\mathcal{L}}}
\newcommand{\Mcal}{{\mathcal{M}}}

\newcommand{\Ucal}{{\mathcal{U}}}

\newcommand{\Acal}{{\mathcal{A}}}

\DeclareMathOperator{\dist}{dist}

\DeclareMathOperator{\supp}{supp}

\DeclareMathOperator{\diverg}{div}

\DeclareMathOperator{\bari}{bar}

\newcounter{fakecounter}

\newcommand{\ind}[1]{\chi_{#1}}

\title{Asymptotic properties of an optimal principal Dirichlet eigenvalue arising in population dynamics}
\author{Lorenzo Ferreri and Gianmaria Verzini}

\begin{document}
\maketitle

\begin{abstract}
We consider a shape optimization problem related to the persistence threshold for a biological 
species, the unknown shape corresponding to the zone of the habitat which is favorable to the 
population. Analytically, this translates in the minimization of a weighted eigenvalue of the 
Dirichlet Laplacian, with respect to a bang-bang indefinite weight. For such problem, we provide 
a full description of the singularly perturbed regime in which the volume of the favorable zone 
vanishes, with particular attention to the interplay between its location and shape.

First, we show that the optimal favorable zone shrinks to a connected, nearly spherical set, 
in $C^{1,1}$ sense, which aims at maximizing its distance from the lethal boundary. Secondly, 
we show that the spherical asymmetry of the optimal favorable zone decays exponentially, 
with respect to a negative power of its volume, in the  $C^{1,\alpha}$ sense, for every 
$\alpha<1$. This latter property is based on sharp quantitative asymmetry estimates for 
the optimization of a weighted eigenvalue problem on the full space, of independent interest.
\end{abstract}%
\noindent
{\footnotesize \textbf{AMS-Subject Classification}}. 
{\footnotesize 49R05, 49Q10, 92D25, 47A75, 35P15
}

\section{Introduction}

Let $\Omega \subset \R^{N}$ denote a bounded domain (open and connected) with regular 
boundary, for instance of class $C^{1,1}$ (even though such assumption can be relaxed, see below). Let $m\in L^\infty(\Omega)$ 
denote a (sign-changing) 
\emph{indefinite weight}. A \emph{principal eigenvalue} of the weighted Dirichlet problem 
\begin{equation}\label{eq:eigprob_intr}
\begin{cases}
-\Delta u = \lambda m u & \text{in }\Omega\\
u=0 & \text{on }\partial\Omega
\end{cases}
\end{equation}
is a real $\lambda$ admitting a positive eigenfunction $u\in H^1_0(\Omega)$. If both the positive and the 
negative parts $m^\pm$ of the weight are nontrivial, \eqref{eq:eigprob_intr} admits two principal eigenvalues, 
and we denote with $\lambda^1(m, \Omega)>0$ the unique positive one:
\begin{equation}\label{eq:eigRay_intr}
\lambda^1(m, \Omega) \coloneqq \inf\left\{\frac{\int_{\Omega}|\nabla u|^2}{\int_{\Omega}m u^2} : 
u \in H^1_0(\Omega),\ \int_{\Omega}m u^2 > 0\right\}
\end{equation}
(in particular, $m\le0$ a.e.\ in $\Omega$ implies $\lambda^1(m, \Omega)=+\infty$).

Among other applications, such eigenvalue is of interest in population dynamics: when modeling the spatial 
dispersal of a population in a heterogenous environment through a reaction-diffusion equation of logistic 
type \cite{Skellam:DisrpersalPopulations,CantrellCosner:indefiniteWeight,Berestycki:PeriodicallyFragmented,MR2191264}, one is lead to consider the 
following evolutive problem for the population density $u=u(x,t)$:
\begin{equation}\label{eq:evol_intr}
\begin{cases}
u_t-d\Delta u = m u - u^2 & x\in\Omega,\ t>0,\\
u = u_0\ge 0 & x\in\Omega,\ t=0,\\
u = 0 & x\in\partial\Omega,\ t>0.
\end{cases}
\end{equation}
Here $m$ describes the heterogeneous habitat, whose favorable and hostile zones for the species 
correspond to positivity and negativity sets of $m$, respectively, and Dirichlet conditions encode a lethal boundary. It is well known that solutions to \eqref{eq:evol_intr} persist (i.e.\ they do not extinguish) if and 
only if $d\lambda^1(m, \Omega)<1$; accordingly, the smaller $\lambda^1(m, \Omega)$ is, the more chances of survival the population has. For this reason, the minimization of $\lambda^1(m, \Omega)$ with respect to the weight, or to other relevant parameters of the model, has been widely investigated. The literature is quite extensive, and we refer to the recent papers \cite{Lamboley:OptimizersRobin,MR3498523,MR3771424,MAZARI2022401,dipierro2021nonlocal} and references therein, both for a more detailed description of the problem, also with different boundary conditions, and for various interesting phenomena, ranging from fragmentation to different nonlocal effects. We also mention \cite{MR1796024,MR2421158}, where related problems arising from the study of composite membranes are investigated.

In this paper we deal with the so-called 
\emph{optimal design problem for the survival threshold}, i.e.\ 
 the minimization of $\lambda^1(m, \Omega)$ with respect to the weight $m$.
For fixed positive constants $\underline{m}$, $\overline{m}$, and $0<\eps<\Lcal(\Omega)$ (where $\Lcal$ denotes 
the Lebesgue measure) let us consider the class of indefinite weights:
\begin{equation}\label{eqn:WeightsClassIndef}
\mathcal{M}_{\varepsilon} \coloneqq \left\{ m \in L^{\infty}(\Omega): -\underline{m} \le m \le \overline{m},\, \int_{\Omega}m = M_{\varepsilon}\right\},
\end{equation}
where $M_{\varepsilon} = \overline{m}\varepsilon - (\mathcal{L}(\Omega)-\varepsilon)\underline{m}$ 
(see below), and the 
class of \emph{bang-bang weights}:
\[
\mathcal{BB}_{\varepsilon} \coloneqq  \left\{ m \in  L^{\infty}(\Omega) : m = \overline{m} \chi_{E}-\underline{m}\chi_{\Omega\setminus E}, \ E \subset \Omega,\ \mathcal{L}(E) = \varepsilon \right\},
\]
where $\ind{A}$ denotes the indicatrix function of the set $A$. 
By definition of $M_\eps$, it is a direct check that $\mathcal{BB}_{\varepsilon}\subset
\mathcal{M}_{\varepsilon}$. Moreover, in the seminal paper \cite{CantrellCosner:indefiniteWeight},
Cantrell and Cosner proved that the 
minimization of $\lambda^1(m, \Omega)$ in the wider class is achieved by an element of the more narrow one: 
\begin{equation}\label{eq:od}
\lambda_{\varepsilon}(\Omega) \coloneqq \inf_{m \in \mathcal{M}_{\varepsilon}} \lambda^1(m, \Omega)  =
\inf_{m \in \mathcal{BB}_{\varepsilon}} \lambda^1(m, \Omega) = \inf_{E \subset \Omega,\ \mathcal{L}(E) = \varepsilon} \lambda^1(E, \Omega), 
\end{equation}
where, with some abuse of notation, for weights in $\mathcal{BB}_{\varepsilon}$ we write 
\[
\lambda^1(E, \Omega)\coloneqq \lambda^1(\overline{m} \chi_{E}-\underline{m}\chi_{\Omega\setminus E}, \Omega).
\]
In particular, the optimal design problem for the survival threshold is turned into a \emph{shape optimization 
problem} involving the favorable zone $E$. In the following, when no confusion arises, referring to the set $\Omega$, we will simply write 
\[
\lambda^1(E)=\lambda^1(E, \Omega),\qquad \lambda_{\varepsilon}=\lambda_{\varepsilon}(\Omega).
\]
According to \cite{CantrellCosner:indefiniteWeight}, it follows that $\lambda_{\varepsilon}$ is 
achieved by a set $E_\eps$, and the first order optimality condition translates in the following 
crucial fact: any optimal set $E_\eps$ is the \emph{super-level set} of the associated positive 
eigenfunction $u_\eps\in H^1_0(\Omega)$. In particular, $E_\eps$ contains the maximum points of $u_\eps$. 

Both for modeling reasons, and from the mathematical point of view, natural questions arise about the \emph{location} and the \emph{shape} of the optimal shape $E_\eps$. This issue is mostly open, and it is completely understood only when the domain $\Omega$ is a ball: in such 
case, by symmetrization arguments it is proved in \cite{CantrellCosner:indefiniteWeight} that $E_\eps$ is a 
concentric ball of measure $\eps$, for every $0<\eps<\Lcal(\Omega)$. On the converse, by 
\cite{Lamboley:OptimizersRobin} if $E_\eps$ is a ball then $\Omega$ itself must be a concentric 
(larger) ball.

In this paper we address the above questions, in the singularly perturbed, small volume regime $\eps\to0$. 
Our aim is to show that, for $\eps$ small, the optimal favorable region $E_\eps$ shrinks to a 
connected, asymptotically spherical shape, which locates at maximal distance from the boundary. 
Moreover, we are going to provide sharp estimates on the optimal eigenvalue 
$\lambda_\eps$, which will reflect on quantitative properties of the shape of the optimal set 
$E_\eps$. On the one hand, the investigation of concentration phenomena to detect the limit 
position and the qualitative asymptotic sphericity of the optimal shape is a classical 
theme in the analysis of variational problems, see e.g. the book \cite{zbMATH01329686} and 
references therein. On the other hand, the use of quantitative estimates to provide explicit 
bounds to the spherical asymmetry of the optimal shape is more recent: we refer for instance 
to \cite{MR2807136}, where this kind of results were obtained for the shape of liquid drops and 
crystals, in the small mass regime, also in the anisotropic setting (where a suitable Wulff 
shape replaces the sphere). It is worth mentioning that, in our case, no perimeter term is present: the main feature of 
the model consists in the interplay between the location and the geometry of the optimal 
favorable set, as the quantitative estimates about the optimal shape are triggered by a 
repulsion effect of the lethal boundary of the box. Incidentally, we also mention that a
first investigation of an anisotropic version of our problem is contained in the recent 
paper \cite{Pellacci:2023aa}.

A crucial role in our description is played by the rescaled limit problem on $\R^N$, which has 
already been studied in \cite[Section 2]{Verzini:Neumann} (see also 
\cite[Proposition 2.1]{FerreriVerzini:AsymptSphericalIndefinite}): consistently  
with the previous notation, for $\Acal\subset\R^N$ (measurable and) bounded let us define the 
principal eigenvalue
\[
\lambda^1(\Acal, \R^N)\coloneqq \inf\left\{\frac{\int_{\R^N}|\nabla v|^2}{\overline{m}\int_{\Acal}v^2 - \underline{m}\int_{\R^N\setminus\Acal} v^2} : 
u \in H^1(\R^N),\ \overline{m}\int_{\Acal}v^2 > \underline{m}\int_{\R^N\setminus\Acal} v^2\right\}.
\]
Then, using symmetrization techniques, one can show that the minimization of such eigenvalue, 
among sets $\Acal$ having prescribed measure, is achieved by a ball, with an associated radial and 
radially decreasing eigenfunction. Precisely, we have that
\begin{equation}\label{eq:lambdazero}
\tilde \lambda_0 \coloneqq \inf_{\Lcal(\Acal)=1} \lambda^1(\Acal, \R^N) = \lambda^1(B, \R^N),
\end{equation}
where $B\subset\R^N$ denotes the ball of unitary Lebesgue measure, centered at the origin, 
and $r_0$ denotes its radius:
\begin{equation}\label{eq:unitaryball}
B=B_{r_0},\qquad \Lcal(B_{r_0}) =1.
\end{equation}
Moreover such minimizer is unique up to translations, and in turn $\lambda^1(B, \R^N)$ is 
achieved by $w\in H^1(\R^N)$, solution of 
\begin{equation}\label{eq:lim_prob}
-\Delta w = \tilde{\lambda}_0 \tilde{m}_0 w \quad \text{in } \R^N ,
\qquad\text{where }\tilde m_0 \coloneqq \overline{m} \chi_{B}-\underline{m}\chi_{\R^N\setminus B};
\end{equation}
here $w>0$, radially symmetric, radially decreasing and normalized in $L^2(\R^N)$, is uniquely 
determined. As a matter of fact, $w$ is explicit in terms of Bessel functions, and it decays 
exponentially at infinity:
\begin{equation}\label{eq:decay_w}
w(\bx) \sim C |\bx|^{-(N-1)/2} e^{-\sqrt{\tilde\lambda_0\underline{m}}|\bx|}\qquad
\text{as }|\bx|\to+\infty.
\end{equation}

Our main results can be divided in two parts. In the first part of the paper 
we analyze the qualitative properties of the optimal favorable set $E_\eps$, as well 
as the asymptotics of the optimal eigenvalue $\lambda_\eps$, leaning on a blow-up procedure. 
In this respect, our main results can be summarized as follows.
\begin{theorem}\label{thm:intro_qual}
There exists $\bar\eps$ such that, for every $0<\eps<\bar\eps$, for every 
optimal eigenfunction $u_\eps>0$ and favorable set $E_\eps$, with $\Lcal(E_\eps)=\eps$, associated to 
$\lambda_\eps=\lambda_\eps(\Omega)=\lambda^1(E_\eps,\Omega)$:
\begin{enumerate}
\item\label{i:unique_max} $u_\eps$ has a unique local maximum point $\bx_\eps$, and in 
particular $E_\eps\ni\bx_\eps$ is connected;
\item\label{i:varphi} $E_\eps$ is radially diffemorphic to a ball centered at $\bx_\eps$, i.e. there exists a function $\varphi_\eps\in C^{1,1}(\sphere^{N-1})$ such 
that ($\Lcal(B_{r_0})=1$)
\[
E_\eps=\left\{\bx: |\bx - \bx_\eps|< \eps^{1/N}\left(r_0+\varphi_\eps\left(
\frac{\bx - \bx_\eps}{|\bx - \bx_\eps|}\right)\right)\right\}.
\]
\setcounter{fakecounter}{\theenumi}
\end{enumerate}
Moreover, as $\eps\to0^+$,
\begin{enumerate}\addtocounter{enumi}{\value{fakecounter}}
\item\label{i:sharp_pos} $\dist(\bx_\eps,\partial\Omega)\to\max_{\bp\in\Omega}\dist(\bp,\partial\Omega)
=:d^*$;
\item\label{i:eig} $\lambda_\eps=\eps^{-2/N}\left(\tilde \lambda_0 +  e^{-2\sqrt{\tilde\lambda_0\underline{m}}(1+o(1))d^*\cdot\eps^{-1/N}}\right)$;
\item\label{i:asimpt_spher} $\|\varphi_\eps\|_{C^{1,1}(\sphere^{N-1})}\to 0$.
\end{enumerate}
Finally, properties \ref{i:varphi}, \ref{i:sharp_pos} and \ref{i:asimpt_spher} remain true if the maximum 
point $\bx_\eps$ is replaced with the barycenter $\bari(E_\eps)$ of the optimal set 
$E_\eps$, and $\varphi_\eps$ with the corresponding polar parametrization of $\partial E_\eps$ 
centered at $\bari(E_\eps)$.   
\end{theorem}

In particular, we obtain an answer to the natural questions raised above: when the volume of 
the favorable zone is small, then such zone is connected, asymptotically spherical, and it 
concentrates at a point whose location tends to maximize its distance to the lethal boundary. 

Notice that, in general, the uniqueness of the pair $u_\eps,E_\eps$ is not guaranteed: 
to disprove it, in view of the theorem, it is sufficient to consider a domain $\Omega$, 
symmetric with respect to a hyperplane, such that points achieving the inradius do not 
belong to such hyperplane. In other words, the properties in Theorem 
\ref{thm:intro_qual}  hold true for every choice of $u_\eps,E_\eps$. Moreover, it will be 
clear from the proof that, using the monotonicity induced by domain inclusion, the regularity 
assumptions on $\partial\Omega$ can be substantially weakened: it suffices to assume that 
$\Omega$ can be approximated, both from inside and from outside, with smooth domains having 
arbitrarily close inradii.

As we mentioned, the proof of Theorem \ref{thm:intro_qual} is based on a blow-up procedure, and 
it exploits 
some analogies with the study of singularly perturbed semilinear elliptic equations. Indeed, 
it is well known in the literature (see \cite{NiTakagi:ShapeLeastEnergy91,NiTakagi:NeumannLocatingPeaks,NiWei:Spike95,MR1736974,BerestyckiWei:SingularlyPerturbedRobin}) that least action solutions of equations like
\[
-\eps^2\Delta u + u = f(u)\qquad\text{in }\Omega,
\] 
complemented with suitable boundary conditions, exhibit concentration phenomena as $\eps\to0$; 
moreover, a sharp asymptotic expansion of the least action level $c_\eps$ allows to locate the 
points where concentration may happen. Now, with respect to such problem, the optimal design problem for the survival threshold 
in population dynamics exhibits several different aspects: first of all, the driving parameter is 
not explicit inside the equation, but it is the prescribed volume of the optimal favorable set; 
moreover, the underlying equation is linear, with a discontinuous and non-homogeneous weight, and 
degenerate on the relevant solutions (the eigenfunctions). On the other hand, the two problems 
present also relevant analogies, triggered by the fact that the unknown shape of 
the optimal favorable set is a superlevel set of the associated eigenfunction, and thus the 
vanishing-volume regime suggests that some concentration for the eigenfunction may happen. Such 
analogies have already been pointed out and exploited in 
\cite{Verzini:AsymptSphericalShapes,Verzini:Neumann}, for the corresponding 
problems with Neumann boundary conditions. 

More specifically, properties \ref{i:unique_max}, \ref{i:sharp_pos} and \ref{i:eig} in Theorem \ref{thm:intro_qual} 
follow by a repeated application of a blow-up argument. In particular, we obtain the uniqueness of 
the maximum points, which reflects on the connectedness of $E_\eps$, by a nontrivial adaptation of 
the arguments in \cite{NiWei:Spike95}; on the other hand, we infer the location of such maxima and 
the sharp expansion of the eigenvalue exploiting a simplified argument inspired by  
\cite{MR1736974}. Incidentally, the blow-up procedure provides information not only on the 
asymptotic behavior of $E_\eps$ and $\lambda_\eps$, but also on that of the optimal eigenfunction 
$u_\eps$. In particular, we show that a suitable rescaling of $u_\eps$ converges, strongly in 
$H^1(\R^N)$ and in $C^{1,\alpha}_\loc$, $\alpha<1$, to the eigenfunction $w$ of the limit problem 
\eqref{eq:lim_prob} (see Lemmas \ref{lem:nohalfspace}, \ref{lem:H1strong}, 
\ref{prop:PropLimitToProblemOverRnIndef} and Corollary \ref{coro:H1strong} ahead for more details).

Concerning the properties of the polar parameterization $\varphi_\eps$ of $\partial E_\eps$, 
it is worth noticing that, being the weight $m_\eps$ merely $L^\infty$, the blow-up procedure 
naturally provides just a $C^{1,\alpha}$-regularity, for $\alpha$ strictly smaller than $1$. To 
recover the $C^{1,1}$ regularity obtained in Theorem \ref{thm:intro_qual}, which is optimal at the 
level of the eigenfunction $u_\eps$, we observe that the partial derivatives of the optimal 
eigenfunction $u_\eps$ satisfy, after a suitable change of variable, a system which can be 
seen as a transmission problem; this crucial remark allows to apply recent regularity results for 
such kind of problems, by Caffarelli, Soria-Carro and Stinga 
\cite{Caffarelli2021:TransmissionProblems} and Dong \cite{Dong2021:TransmissionProblems}, thus 
concluding the proof of Theorem \ref{thm:intro_qual}. We remark that such further regularity is not 
only interesting by itself, but it also allows to simplify some arguments that we use in the second 
part of the paper.\medskip

Properties \ref{i:varphi}, \ref{i:asimpt_spher} in Theorem \ref{thm:intro_qual} assert that, from a qualitative 
point of view, the shape of the optimal favorable set $E_\eps$ is asymptotically a shrinking sphere, in a $C^{1,1}$ 
sense. Then a natural question, that we address in the second part of the paper, is to which extent such assertion can be made quantitative. The answer is contained in the following result.
\begin{theorem}\label{thm:intro_quant_domain}
Under the assumptions and notation of Theorem \ref{thm:intro_qual}, let $0<\eps<\bar\eps$
and let $\varphi_\eps$ denote the polar parametrization of  $\partial E_\eps$, as in 
property \ref{i:varphi} of such theorem, centered at $\bari(E_\eps)$. 

Then, for a possibly smaller value of $\bar\eps$ and for every $0<\alpha<1$ there exist 
positive constants $M,C$ such 
that 
\[
\| \varphi_{\varepsilon} \|_{C^{1,\alpha}(\sphere^{N-1})} \le C e^{-M\eps^{-1/N}},
\]
for every $0<\eps<\bar \eps$.
\end{theorem}

The proof of Theorem \ref{thm:intro_quant_domain} is based on a sharp 
quantitative form of the inequality provided by the limit problem \eqref{eq:lambdazero}, 
\eqref{eq:lim_prob}, which is of independent interest. Such result is expressed 
in terms of a quantification of the spherical asymmetry of the favorable 
set, based on the following definition.
\begin{definition}\label{def:nearly_spherical}
A (bounded) set $\Acal\subset\R^N$ is \emph{nearly spherical of class $C^{1,\alpha}$ parametrized by $\varphi$}, 
if there exists $\varphi\in C^{1,\alpha}(\sphere^{N-1})$, with $\|\varphi\|_{L^\infty}\le r_0/2$, such 
that $\partial \Acal$ is represented as
\[
 \partial \Acal =\left\{\bx\in\R^N:\bx = (r_0 + \varphi(\theta))\theta\text{ for }\theta\in\sphere^{N-1}
\right\}
\]
(recall that $r_0$ is the radius of the ball of unit measure in $\R^N$: $B=B_{r_0}(\mathbf{0})$).
\end{definition}

Under the above definition, and denoting with $\bari(\mathcal{A})$ the barycenter of any bounded 
measurable set in $\R^N$, we obtain the following sharp quantitative version of 
problem \eqref{eq:lambdazero}.
\begin{theorem}\label{thm:quantitStabRNSharpPos}
There exist constants $C>0$, $\delta>0$ such that, for any $C^{1, 1}$ nearly spherical set $\mathcal{A} \subset \R^N$ satisfying
\begin{itemize}
    \item[(i)] $\bari(\mathcal{A}) = \mathbf{0}$, 
    \item[(ii)] $\mathcal{L}(\mathcal{A}) = 1$
    \item[(iii)] $\| \varphi_{\mathcal{A}} \|_{C^{1, 1}(\sphere^{N-1})} \le \delta$,
\end{itemize}
it holds
\[
\lambda^1(\mathcal{A}, \R^N) - \lambda^1(B, \R^N) \ge C \| \varphi_{\mathcal{A}} \|^2_{L^2(\sphere^{N-1})} .
\]
\end{theorem}
In the above theorem, the constant $C>0$ is obtained by a constructive method and in 
principle it is computable. As we mentioned, the above estimate is sharp, in the sense that we 
prove also the existence of a larger constant $C'>C$ such that, under the same assumptions, 
also the reverse inequality holds:
\begin{equation}\label{eq:sharp_indeed}
\lambda^1(\mathcal{A}, \R^N) - \lambda^1(B, \R^N) \le C' \| \varphi_{\mathcal{A}} \|^2_{L^2(\sphere^{N-1})} 
\end{equation}
(see Lemma \ref{lem:sharp_indeed} ahead). We also remark that the constraint on the barycenter of the nearly spherical set is necessary to obtain 
that $B$ is locally a \emph{strict} minimum shape, or equivalently that $C>0$ strictly. Indeed, since 
problem \eqref{eq:lambdazero} is settled in the whole $\R^N$, translations of $\partial B$ 
are still global minimizers. 

Theorem \ref{thm:quantitStabRNSharpPos} can be easily improved, at the regularity level, using the 
Gagliardo-Nirenberg inequality on the sphere $\sphere^{N-1}$ (and in turn, such improvement is 
crucial to deduce Theorem \ref{thm:intro_quant_domain} from property  \ref{i:eig} of Theorem 
\ref{thm:intro_qual}).
\begin{corollary}\label{coro:quantitative_GN}
Under the assumptions of Theorem \ref{thm:quantitStabRNSharpPos}, for every 
$0<\alpha<1$ there exists a constant $C>0$, depending also on $\delta$, such that  
\[
\lambda^1(\mathcal{A}, \R^N) - \lambda^1(B, \R^N) \ge C 
\| \varphi_{\mathcal{A}} \|_{C^{1,\alpha}(\sphere^{N-1})}^{(4+N)/(1-\alpha)} .
\]
\end{corollary}

The proof of Theorem \ref{thm:quantitStabRNSharpPos} is carried out 
expanding the map 
\[
\mathcal{A}\mapsto\lambda^1(\mathcal{A}, \R^N)
\] 
in Taylor series, centered at $B$, up to the second order. After the seminal paper 
\cite{Fuglede_1989}, this is a by now classical argument 
in shape optimization problems, see for instance \cite{DambrinePierre:StabilityEquilibriumShapes,Dambrine2002:ShapeHessian,HenrotPierre:ShapeOptimiz,   Dambrine2011:shapeSensitivity}, also to deduce sharp quantitative estimates for volumetric 
geometric-functional inequalities as in \cite{Brasco:SharpFaberKrahn, DePhilippis:SharpIsocapacitary}. 
In particular, we advise the paper \cite{Mazari2020:QuantitaiveShrodinger}, where Mazari considers 
the optimization of the first eigenvalue of a Schr\"odinger operator in the ball, with respect to 
a bang-bang potential, obtaining results related with ours. In particular, as remarked in 
\cite[Sec. 2.5.2]{Mazari2020:QuantitaiveShrodinger}, in shape optimization problems where the shape 
is a subdomain of a fixed domain, the natural norm to prove quantitative asymmetry estimates for nearly spherical sets 
seems to be the $L^2(\sphere^{N-1})$ norm of the deformation, rather than the 
$H^{1/2}(\sphere^{N-1})$ norm as it is more usual in the literature 
(see e.g.\  \cite{DambrineLamboley:stabShapeOptim, Brasco:SharpFaberKrahn, DePhilippis:SharpIsocapacitary, Dambrine2002:ShapeHessian}). In this respect, our 
results corroborate this observation. 
On the other hand, with respect to \cite{Mazari2020:QuantitaiveShrodinger}, in our proof 
we can take advantage of the $C^{1,1}$ regularity of the shape to simplify 
part of the argument. In particular, on the lines of \cite{Mazari2020:QuantitaiveShrodinger}, 
we think that Theorem \ref{thm:quantitStabRNSharpPos} may be proved also starting from a 
weaker $C^{1,\alpha}$ regularity assumption, even though we do not pursue this issue here.

The paper is structured as follows:

\tableofcontents

\subsection*{Notation}

\begin{itemize}
\item $\Lcal(A)$ denotes the $N$-dimensional Lebesgue measure of the set $A\subset\R^N$.  
\item $B\subset \R^N$ denotes the ball centered at the origin, with $\Lcal(B)=1$; consequently, for every $t>0$, $\Lcal(tB)=t^N$. 
Moreover, $r_0$ is defined in such a way that
\[
B=B_{r_0}(\mathbf{0}).
\]
\item $C$, $C'$ and so on denote positive constants, which may change the value from line to line.
\end{itemize}

\section{Blow-up procedure}\label{section:blowUpIndef}

In this section we begin to study the qualitative asympotics for  $\lambda_{\varepsilon}$, 
the first generalized eigenvalue with indefinite weight. We proceed by means of a blow-up 
argument, similarly to \cite[Section 4]{Verzini:Neumann} and 
\cite[Section 3]{FerreriVerzini:AsymptSphericalIndefinite},  to which we refer for further details.  

As recalled in the introduction, we denote with $u_{\varepsilon}$ a family of 
positive eigenfunctions associated to $\lambda_{\varepsilon}$, normalized in $L^2(\Omega)$, having 
associated optimal weight $m_\eps\in L^\infty(\Omega)$ and optimal favorable $E_\eps\subset\Omega$, 
with $\Lcal(E_\eps)=\eps$. 
Although for our starting purposes the blow-up procedure will be centered at 
$\mathbf{x}_{\varepsilon} \in \Omega$, a global maximum point of $u_{\varepsilon}$, 
later we will need to choose different centers. For this reason it is convenient to 
illustrate the procedure for an arbitrary choice of the blow-up centers $\bp_\eps$.

Let 
\begin{equation}\label{eq:kbeta}
\bp_\eps\in \Omega,\qquad
k_{\varepsilon} \coloneqq \varepsilon^{1/N}.
\end{equation}
We define the blow-up functions
\begin{equation}\label{eq:BlowUpFuncs}
\tilde{u}_{\varepsilon}(\mathbf{x}) \coloneqq k_{\varepsilon}^{N/2} \, u_{\varepsilon}
(\mathbf{p}_{\varepsilon} + k_{\varepsilon} \mathbf{x}), \qquad
\tilde{m}_{\varepsilon}(\mathbf{x}) \coloneqq m_{\varepsilon}(\mathbf{p}_{\varepsilon} + k_{\varepsilon} \mathbf{x}),
\end{equation}
and (superlevel) sets
\begin{equation}\label{eq:BlowUpSets}
\begin{split}
\tilde{\Omega}_{\varepsilon} &\coloneqq \left\{ \mathbf{x} \in \R^N : \mathbf{p}_{\varepsilon} + k_{\varepsilon} \mathbf{x} \in \Omega \right\}=\left\{ \mathbf{x} \in \R^N : \tilde{u}_{\varepsilon}(\mathbf{x})>0 \right\},\\
\tilde E_{\varepsilon} &\coloneqq \left\{ \mathbf{x} \in \R^N : \mathbf{p}_{\varepsilon} + k_{\varepsilon} \mathbf{x} \in E_{\varepsilon} \right\}=
\left\{ \mathbf{x} \in \R^N : \tilde{u}_{\varepsilon}(\mathbf{x})>\alpha_\eps \right\},
\end{split}
\end{equation}
where $\alpha_\eps= k_{\varepsilon}^{N/2}\left. u_{\varepsilon}\right|_{\partial E_\eps}>0$. The normalizations are chosen in such a way that, for every $\eps>0$,
\[
\int_{\tilde\Omega_\eps} \tilde u_\eps^2 =\mathcal{L}\left( \tilde E_{\varepsilon} \right) = 1.
\]
Moreover, the functions $\tilde{u}_{\varepsilon}$ solve, in $H_0^1(\tilde{\Omega}_{\varepsilon})$,
\begin{equation}\label{eq:GeneralDiffProblemBUIndef}
\begin{cases}
-\Delta \tilde{u}_{\varepsilon} = \tilde{\lambda}_{\varepsilon} \tilde{m}_{\varepsilon} \tilde{u}_{\varepsilon} & \text{in } \tilde{\Omega}_{\varepsilon} , \\
\tilde{u}_{\varepsilon} = 0 & \text{on } \partial\tilde{\Omega}_{\varepsilon} ,
\end{cases}
\end{equation}
where 
\begin{equation}\label{eqn:DefBlowUpEigenvIndef}
    \tilde{\lambda}_{\varepsilon}=\tilde{\lambda}_{\varepsilon}(\tilde{\Omega}_{\varepsilon}) \coloneqq 
    k_{\varepsilon}^2\lambda_{\varepsilon}(\Omega)
\end{equation}
satisfies
\[
\tilde{\lambda}_{\eps} = \inf_{A \subset \tilde{\Omega}_\eps : \mathcal{L}(A) = 1} \lambda^1(A, \tilde{\Omega}_\eps)= \lambda^1(\tilde E_\eps, \tilde{\Omega}_\eps).
\]

Since $\Omega$ is bounded, from any sequence $\varepsilon_n \to 0$ we can extract a subsequence, still 
denoted by $\varepsilon_n$, such that $\mathbf{p}_{\varepsilon_n} \to \mathbf{p}_\infty \in 
\overline{\Omega}$ for $n \to +\infty$. For easier notation we write 
$\tilde{\lambda}_n = \tilde{\lambda}_{\varepsilon_n}$, $\tilde{u}_n = \tilde{u}_{\varepsilon_n}$, 
$\tilde{m}_n = \tilde{m}_{\varepsilon_n}$, and so on. In the following we are going to 
prove convergence of such sequences, possibly up to further subsequences. The blow-up limits will be 
related to the corresponding counterpart of the limit problem \eqref{eq:lambdazero}, 
\eqref{eq:lim_prob}, namely 
$\tilde\lambda_0$, $w$, $\tilde m_0$. We remark that, in case the limit of the subsequence is  
independent of the starting sequence, then convergence of the whole family follows, as $\eps\to0$.

We start identifying the limit of the blow-up sequences. To this aim, we notice that, by trivial 
extension, we can assume that $\tilde u_n\in H^1(\R^N)$. On the other hand, it is convenient to 
extend $\tilde m_n$ outside $\tilde \Omega_n$ as $-\underline{m}$. In this way we obtain that 
$\tilde m_n+\underline{m} = (\overline{m}+\underline{m})\chi_{\tilde E_n}$ so that
\begin{equation}\label{eq:M'}
\tilde m_n \in \Mcal':=\left\{m\in L^\infty(\R^N):-\underline{m} \le  m  \le \overline{m} , \quad \int_{\R^N}( m +\underline{m}) \le (\underline{m} + \overline{m}) \right\}.
\end{equation}
\begin{lemma}\label{lem:conv2la_0Indef}
For any choice of the blow-up centers $(\bp_n)_n$ we have, as $n \to +\infty$,
\begin{enumerate}
\item $\tilde{\lambda}_n \to \tilde{\lambda}_0$;
\item up to subsequences, 
\[
\tilde\Omega_n\to\Omega_\infty=
\begin{cases}
\R^N & \text{if }d(\mathbf{p}_n, \partial\Omega)/k_n \to +\infty\\
\text{an open half-space }H  & \text{if }d(\mathbf{p}_n, \partial\Omega)/k_n \le C,
\end{cases}
\]
in the sense that, for every $K\subset\subset H$ (resp. $\R^N\setminus H$), we have $K\subset\subset\tilde\Omega_n$ (resp. 
$\R^N\setminus \tilde\Omega_n$) for $n$ sufficiently large.
\item There exists $ m_\infty \in \Mcal'$ such that, up to subsequences, $\tilde{m}_n\to m_\infty$ weakly-$\ast$ in $L^{\infty}(\R^N)$.
\item There exists $ u_\infty  \in H^1(\R^N)\cap H^1_0(\Omega_\infty) $ such that, up to subsequences, 
$\tilde u_n\to u_\infty$ both weakly in $H^1(\R^N) $ and $C^{1, \alpha}_{\loc}(\Omega_\infty)$, 
for every $0 < \alpha < 1$, and
\begin{equation}\label{eq:infty}
\int_{\R^N} \nabla  u_\infty  \cdot\nabla v = \tilde{\lambda}_0 \int_{\R^N}  m_\infty   u_\infty  v \qquad \forall v \in H^1_0(\Omega_\infty).
\end{equation}
Moreover, $u\ge0$ in $\R^N$.
\end{enumerate}
\end{lemma}
\begin{proof}
The proof of 1 is analogous to that of \cite[Lemma 3.1]{FerreriVerzini:AsymptSphericalIndefinite} 
(notice, in particular, that $\tilde{\lambda}_n\ge\tilde{\lambda}_0$ is bounded above by the corresponding generalized 
eigenvalue with spherical weight). 

Claim 2 follows by regularity of $\partial\Omega$: 
in particular, $H$ depends in an elementary way on $\bp_\infty$, on the tangent plane to 
$\partial\Omega$ at $\bp_\infty$ and on the asymptotic direction 
$\lim_n\frac{\bp_n-\bp_\infty}{|\bp_n-\bp_\infty|}$.

To show 3 it is sufficient to notice that $(\tilde m_n)_n$ is bounded in $L^\infty(\R^N)$, and 
that $\Mcal'$, defined in \eqref{eq:M'}, is closed with respect to the weak-$*$ topology (see 
\cite[Lemma 4.1, Step 3]{Verzini:Neumann} for further details). 

To prove 4 one can proceed similarly as in \cite[Proposition 3.3]{FerreriVerzini:AsymptSphericalIndefinite}: first, using the equation and point 1, we show 
the boundedness of $(\tilde u_n)_n$, obtaining weak $H^1(\R^N)$, strong $L^2_\loc(\R^N)$ and a.e. 
convergence to $u_\infty\in H^1(\R^N)$, $u_\infty\ge0$ (recall that $\tilde u_n$ is positive and 
$L^2$-normalized in $\tilde\Omega_n$). Next, by point 2 and definition of $\tilde u_n$, we have 
that $u_\infty = 0$ a.e. in any $K\subset\subset \R^N\setminus \overline{\Omega}_\infty$, implying that (the 
restriction to $\Omega_\infty$ of) $u_\infty$ belongs to $H^1_0(\Omega_\infty)$. On the other hand, 
choose an arbitrary $K\subset\subset \Omega_\infty$ and take any $v \in C_0^{\infty}(K)$. 
Then, thanks to 2, 
for $n$ large enough $\supp(v) \subset \tilde{\Omega}_n$. Thus $\tilde{u}_n$ satisfies:
\[
\int_{\R^N} \nabla \tilde{u}_n \cdot \nabla v = \tilde{\lambda}_n \int_{\R^N} \tilde{m}_n \tilde{u}_n v  .
\]
Using 1, 3 and the density of $C_0^{\infty}(\Omega_\infty)$ in $H^1_0(\Omega_\infty)$ we infer 
\eqref{eq:infty}. Finally, compactness in $C^{1, \alpha}(K)$, $0<\alpha<1$, follows by the 
weak equation above, using elliptic regularity and a standard bootstrap argument (see 
e.g. \cite[Lemma 4.1, Step 1]{Verzini:Neumann}). 
\end{proof}

In the following we focus on a subsequence, still denoted with $\eps_n$, such that all the claims 
in Lemma \ref{lem:conv2la_0Indef} hold true.
\begin{lemma}\label{lem:nohalfspace}
For any choice of the blow-up centers $(\bp_n)_n$ we have that one of the following two 
alternatives hold:
\begin{enumerate}
\item either $u_\infty\equiv0$ in $\R^N$,
\item or $\Omega_\infty=\R^N$ and, up to (the same) translation,
\[
m_\infty = \tilde m_0,\qquad u_\infty  = A w,
\]  
for some positive constant $A$. In particular, $\dist(\bp_n,\partial\Omega)/k^{1/N}\to+\infty$ as $n\to+\infty$.
\end{enumerate}
\end{lemma}
\begin{proof}
Assume that $u_\infty$ is not the null function. Then, taking $v=u_\infty$ in \eqref{eq:infty} 
we obtain
\[
0<\int_{\R^N} |\nabla  u_\infty|^2 = \tilde{\lambda}_0 \int_{\R^N}  m_\infty   u_\infty^2,
\qquad\text{i.e. }
\tilde{\lambda}_0 = \frac{\int_{\R^N}|\nabla  u_\infty |^2}{\int_{\R^N}  m_\infty   u_\infty ^2} .
\]
Since $u\in H^1(\R^N)$, $m\in\Mcal'$ and $\int_{\R^N}  m_\infty   u_\infty ^2>0$, the lemma follows 
by the variational characterization of $\tilde\lambda_0$, see \cite[Thm. 2.2]{Verzini:Neumann}.
\end{proof}

To conclude we recall that, in view of \cite{FerreriVerzini:AsymptSphericalIndefinite}, one can
obtain strong $H^1$ convergence of the blow-up sequence under the further assumption that the 
optimal favorable sets $\tilde E_n$ are uniformly bounded.
\begin{lemma}\label{lem:H1strong}
Under the above notations, assume that $\tilde u_n (\mathbf{0})\ge C>0$ and that 
$\tilde E_n \subset K$, for a suitable constant $C$ and a suitable compact set $K$. Then
\[
\tilde u_n \to w \qquad\text{strongly in }H^1(\R^N)\text{ and in }C^{1,\alpha}_\loc
\]
up to a translation.
\end{lemma}
\begin{proof}
By Lemmas \ref{lem:conv2la_0Indef} and \ref{lem:nohalfspace} we only need to prove the strong $H^1$ convergence, which follows 
proceeding exactly as in \cite[Section 4]{FerreriVerzini:AsymptSphericalIndefinite}.
\end{proof}

\section{Qualitative properties of the favorable region}\label{section:connectedness}

In this section we prove some qualitative geometric properties of the favorable region $E_\eps$, 
and of its counterpart in the blow-up scale $\tilde E_\eps$, in the asymptotic limit 
$\varepsilon \to 0^{+}$. In particular, we will show that it is connected and that it 
converges in a suitable sense to a ball. In particular we will show property 
\ref{i:unique_max} of Theorem \ref{thm:intro_qual}, as well as a version in $C^{1,\alpha}$, 
$\alpha<1$, of properties \ref{i:varphi} and \ref{i:asimpt_spher}.

The main arguments are based on the refinement of the blow-up analysis introduced in the previous 
section, when centered at maximum points for $u_\eps$. In particular, throughout this section we 
take $\eps_n\to0$,
\[
\bp_n=\bx_n=\bx_{\eps_n}\text{ a global maximum point of }u_{\eps_n},
\]
and we denote with $\tilde u_n$ the corresponding blow-up sequence as in \eqref{eq:BlowUpFuncs}. 
\begin{lemma}\label{prop:PropLimitToProblemOverRnIndef}
Choosing $\bp_n = \bx_n$ we obtain that the second alternative in Lemma \ref{lem:nohalfspace} 
occurs, and no translation is needed.
\end{lemma}
\begin{proof}
Taking into account Lemmas \ref{lem:conv2la_0Indef} and \ref{lem:nohalfspace} we have to 
show that 
\begin{equation}\label{eq:maxbddbelow}
\tilde u_n(\mathbf{0}) \ge C>0.
\end{equation}
Indeed, this will imply convergence of $\tilde u_n$ to $Aw$, $A>0$, up to a translation. However, 
by uniform convergence (point 4 of Lemma \ref{lem:conv2la_0Indef}) and since by construction 
$\tilde{u}_n(\mathbf{0})$ is a maximum, necessarily the translation must be $\mathbf{0}$ (recall 
that $w$ is radially symmetric and decreasing). 

On the other hand, \eqref{eq:maxbddbelow} easily follows by the $L^2$ normalization, as in 
\cite[Proposition 3.3 (iv)]{FerreriVerzini:AsymptSphericalIndefinite}. Indeed, since
\[
0< \overline{m}\int_{\tilde E_n} \tilde u_n^2 - \underline{m}\int_{\R^N\setminus\tilde E_n} \tilde u_n^2  = (\overline{m}+ \underline{m})\int_{\tilde E_n} \tilde u_n^2 - \underline{m},
\]
$\mathbf{0}$ is a global maximum point and $\Lcal(\tilde E_n)=1$, we infer
\[
\tilde u_n^2(\mathbf{0}) \ge \int_{\tilde E_n} \tilde u_n^2 \ge 
\frac{\underline{m}}{\overline{m}+ \underline{m}}>0.
\qedhere
\]
\end{proof}
\begin{remark}\label{rem:change_center_1}
We notice that the above lemma (and all its consequences throughout this section) 
hold true by choosing different centers for the blow-up procedure, say  
$\bp_n=\by_n$, as long as $|\bx_n-\by_n |/k_n\to0$ as $n\to+\infty$.
\end{remark}

The main properties of $\tilde E_n$ will descend by the fact that it is a superlevel set.
Recall that, as defined in \eqref{eq:BlowUpSets}, $\alpha_n$ is the value of $\tilde{u}_n$ at the boundary of the favorable set. Accordingly, let us denote with $\alpha_0$ the positive real value 
such that 
\[
\mathcal{L}\left( \left\{  u_\infty  > \alpha_0 \right\} \right) = 1,
\] 
where $ u_\infty = Aw$ is as in Lemma \ref{lem:nohalfspace}. We have the following key fact.
\begin{lemma}\label{lemma:levelSetValueConvergenceIndef}
We have $\alpha_n \to \alpha_0$ as $n\to+\infty$.
\end{lemma}
\begin{proof}
First notice that  $0<\alpha_n\le\tilde{u}_n(\mathbf{0}) \to u_\infty(\mathbf{0})$. As a consequence, 
we can extract a subsequence, still denoted with $\alpha_n$, such that $\alpha_n \to \bar{\alpha}$ 
for $n \to +\infty$. To conclude, we show by contradiction that both $\bar{\alpha} \ge \alpha_0$ and $\bar{\alpha} \le \alpha_0$.

Suppose $\bar{\alpha} < \alpha_0$. Then, since as already remarked $ u_\infty $ is radially strictly decreasing, we have that $\mathcal{L}(\{  u_\infty  > \bar{\alpha} \}) > 1$. Thus, by local uniform convergence on compact sets, for $n$ large enough it also holds that $\mathcal{L}(\{ \tilde{u}_n > \alpha_n \}) > 1$, which is a contradiction with the normalization  $\mathcal{L}(\{ \tilde{u}_n > \alpha_n \}) =1$. Thus $\bar{\alpha} \ge \alpha_0$.

Now suppose that $\bar{\alpha} > \alpha_0$. Consider the function $\chi_{\{  u_\infty  > \alpha_0 \}} \in L^1(\R^N)$. By Proposition \ref{prop:PropLimitToProblemOverRnIndef} (ii) we have that
\[
\int_{\R^N} \tilde{m}_n \chi_{\{  u_\infty  > \alpha_0 \}} \to \int_{\R^N} \tilde{m}_0 \chi_{\{  u_\infty  > \alpha_0 \}} = \overline{m} \qquad \text{for } n \to +\infty.
\]
On the other hand, by uniform convergence on compact sets we have, for $ n \to +\infty$,
\[
\begin{split}
   \int_{\R^N} \tilde{m}_n \chi_{\{  u_\infty  > \alpha_0 \}} &= \overline{m} \mathcal{L}(\{ \tilde{u}_n > \alpha_n \} \cap \{  u_\infty  > \alpha_0 \}) - \underline{m} \mathcal{L}(\{ \tilde{u}_n \le \alpha_n \} \cap \{  u_\infty  > \alpha_0 \})  \\
    & \le \overline{m} \mathcal{L}(\{ \tilde{u}_n > \alpha_n \} \cap \{  u_\infty  > \alpha_0 \}) \le \overline{m} \mathcal{L}(\{ u_\infty > \bar{\alpha} \}) + o(1) < \overline{m},
\end{split}
\]
which leads to a contradiction. 

Thus $\alpha_n \to \bar{\alpha} = \alpha_0$ for $n \to +\infty$ and the proof is concluded.
\end{proof}

To prove the connectedness of the favorable region for $\eps$ small we proceed by contradiction, 
assuming that $\tilde{u}_n$ admits a second (possibly relative) maximum point, that we denote with 
$\mathbf{x}_{\varepsilon_n}'=\mathbf{x}_{n}'$, such that
\[
\bx'_n \neq \bx_n.
\]
Indeed, this would be the case if we could build a sequence $\varepsilon_n\to0$ such that $E_{n}$ is \emph{not} connected, for any $n \in \N$, with a maximum point in each connected component.

The first step to obtain a contradiction consists in performing the blow-up analysis, centered at 
$\bx'_n$.
\begin{lemma}\label{lem:othermax}
Let us choose $\bp_n = \bx'_n$ and denote with $\tilde u'_n$, $\tilde m'_n$  the corresponding blow up 
sequences. Then, up to subsequences,
\[
\tilde m'_n \to \tilde m_0\text{ weakly-$*$ in }L^\infty(\R^N),\qquad \tilde u'_n  \to A' w\text{ weakly in }H^1(\R^N)\text{ and in }C^{1,\alpha}_\loc,
\]  
$A'>0$, as $n\to+\infty$.
\end{lemma}
\begin{proof}
Again, the proof will follow from Lemmas \ref{lem:conv2la_0Indef} and \ref{lem:nohalfspace} once we 
know that 
\begin{equation*}
\tilde u'_n(\mathbf{0}) \ge C>0.
\end{equation*}
This follows by Lemma \ref{lemma:levelSetValueConvergenceIndef}, since $\bx'_n\in  E_n$ and thus
\[
\tilde u'_n(\mathbf{0})> \alpha_n \to \alpha_0>0.
\qedhere
\]
\end{proof}

Now, only three mutually exclusive scenarios can happen for $n \to +\infty$:
\begin{enumerate}
    \item \emph{(repulsion)} $ \quad \dfrac{\dist(\mathbf{x}_{n}', \bx_n)}{k_n} \to +\infty$,
    \item \emph{(coexistence)} $ \quad c \le \dfrac{\dist(\mathbf{x}_{n}', \bx_n)}{k_n}\le C$, for some positive constants $c, C$ independent of $n$,
     \item \emph{(collapsing)} $ \quad \dfrac{\dist(\mathbf{x}_{n}', \bx_n)}{k_n} \to 0$.
\end{enumerate}
Our aim is to prove that none of the cases above is possible, thus leading to a contradiction.
\begin{lemma}\label{lemma:noRepulsionIndef}
In the limit $n \to +\infty$ repulsion cannot happen.
\end{lemma}
\begin{proof}
We proceed by contradiction. So, suppose repulsion.

By Lemma \ref{lem:othermax} we have that the blow-up sequence $\tilde{u}_n'$ with respect to the points $\mathbf{x}_{n}'$ converges uniformly on compact sets to $u_\infty '=A'w$, while as usual 
$\tilde{u}_n'$ converges to $u_\infty =Aw$. We now prove that actually $ u_\infty  =  u_\infty '$.

Let us denote with $\alpha_0'$ the positive real value such that $
\mathcal{L}\left( \left\{  u_\infty ' > \alpha_0' \right\} \right) = 1$. Then one can repeat 
the proof of Lemma \ref{lemma:levelSetValueConvergenceIndef} to obtain that $\alpha_n \to \alpha_0'$ for $n \to +\infty$. Thus $\alpha_0 = \alpha_0'$ and thanks to their definitions, and to the strict monotonicity of $w$, this implies $ u_\infty ' =  u_\infty $.

Now the proof follows easily. Indeed, using the uniform convergence on compact sets, for any $\eta> 0$, for $n$ large enough there exists an $R > 0$ independent of $n$ such that
\[
\mathcal{L}\left(\{ \tilde{u}_n > \alpha_n \} \cap B_R(\mathbf{0})\right) > 1 - \eta \qquad \text{and} \qquad \mathcal{L}\left(\{ \tilde{u}_n' > \alpha_n \} \cap B_R(\mathbf{0}) \right) > 1 - \eta.
\]
However, due to the repulsion assumption, this implies that 
$\mathcal{L}(\{ \tilde{u}_n > \alpha_n \} ) > 2(1 - \eta)$, which is a contradiction.
\end{proof}
\begin{lemma}\label{lemma:noCoexistenceIndef}
In the limit $n \to +\infty$ coexistence cannot happen.
\end{lemma}
\begin{proof}
We proceed again by contradiction. Notice that, by assumption, 
\[
\bz_n \coloneqq \frac{\bx'_n-\bx_n}{k_n} \subset K \subset\subset\R^N,\qquad \bz_n\to\bz_\infty\neq\mathbf{0},
\]
up to subsequences , as $n\to+\infty$.

Moreover, by definition, $\nabla \tilde u_n(\mathbf{0}) = \nabla \tilde u_n(\mathbf{z}_{n}) = 
\mathbf{0}$. Since $\tilde{u}_n$ converges in $C^1(K)$  to $Aw$, $A>0$ (see Lemma \ref{lem:nohalfspace}), we infer that $\nabla w(\mathbf{0}) =\nabla w(\mathbf{z}_\infty) = \mathbf{0}$ with $\bz_\infty\neq\mathbf{0}$, contradiction. 
\end{proof}
\begin{lemma}\label{lemma:noCollapsingIndef}
In the limit $n \to +\infty$ collapsing cannot happen.
\end{lemma}
\begin{proof}
We are in a similar situation as in the previous lemma, but now
\[
\bz_n \coloneqq \frac{\bx'_n-\bx_n}{k_n}\neq \mathbf{0},\qquad \bz_n\to\mathbf{0},
\]
up to subsequences , as $n\to+\infty$. Then we obtain a contradiction using \cite[Lemma 4.2]{NiTakagi:ShapeLeastEnergy91} with $\phi =  u_\infty $ and $\psi = \tilde{u}_n$ once one notices that, 
in Lemma \ref{lem:conv2la_0Indef}  due to standard elliptic regularity theory, the convergence in $B_a(\mathbf{0})$ for $a > 0$ small enough is actually in $C^k$ for any $k \ge 0$.
\end{proof}
\begin{corollary}\label{coro:1di1}
Property \ref{i:unique_max} of Theorem \ref{thm:intro_qual} holds true.
\end{corollary}

Once we know that, for $n$ sufficiently large, $\tilde u_n$ has a unique local maximum point and 
$\tilde E_n$ is connected, we can deduce that the optimal set is asymptotically spherical, in the following sense. 
\begin{lemma}\label{lemma:favRegionAsymptSphericalIndef}
For every $\eta>0$ sufficiently small there exists $\bar n$ such that, for every $n>\bar n$,
\[
 (1-\eta)B\subset \tilde E_{n}\subset (1+\eta)B.
\]
\end{lemma}
\begin{proof}
It follows readily from Lemmas \ref{lem:H1strong} and  
\ref{lemma:levelSetValueConvergenceIndef} (ii), since they imply that $\partial \tilde E_{n}$ is eventually contained in any neighborhood of $\partial B$ (i.e.\ of the $\alpha_0$ level set of $ u_\infty $).
\end{proof}

Once $\tilde E_n$ is uniformly bounded, Lemma \ref{lem:H1strong} applies straightforwardly. 
\begin{corollary}\label{coro:H1strong}
Under the above notation, $u_\infty = w$ and the convergence is strong in $H^1(\R^N)$.
\end{corollary}
Moreover, exploiting the $C^{1, \alpha}$ regularity of the family 
$\tilde{u}_{n}$, we can prove that $\partial \tilde E_n$ is nearly spherical, in the sense of Definition \ref{def:nearly_spherical}.
\begin{proposition}\label{prop:FavRegC1aRegSharpIndef}
For $n$ sufficiently large $\partial \tilde{E}_{n}$ is nearly spherical of class $C^{1, \alpha}$, for all $0 \le \alpha < 1$, represented by  $\varphi_{\tilde{E}_{n}}$. Moreover
\begin{equation}\label{eqn:NormTo0FavRegC1aRegSharpInde}
\| \varphi_{\tilde{E}_{n}} \|_{C^{1, \alpha}(\partial B)} \to 0 ,\qquad \text{for } n \to +\infty.
\end{equation}
\end{proposition}
\begin{proof}
We use polar coordinates $\bx=\rho\theta$, with $\rho>0$, $\theta\in \sphere^{N-1}$. Recalling that $\partial\tilde E_n$  
is the $\alpha_n$-level of the function $\tilde u_n$, we obtain that 
$\varphi_{\tilde E_n} = \rho - r_0$ is implicitly defined as
\[
\tilde u_\eps((r_0 + \varphi_{\tilde E_n}) \theta) = \alpha_n.
\]
Now, consider $n$ large enough so that $\partial \tilde E_n \subset B_R\setminus B_r$ and 
\begin{equation}\label{eqn:eqn1LevelC1aSurSharpIndef}
\max_{\overline{B_R\setminus B_r}}\partial_{\rho} \tilde{u}_{n} \le \frac{1}{2} \max_{\overline{B_R\setminus B_r}}\partial_{\rho} w < 0,
\end{equation}
where $r<r_0<R$ are fixed. This is possible thanks to Lemmas \ref{lem:H1strong} 
and \ref{lemma:favRegionAsymptSphericalIndef}. Since 
$\partial \tilde E_n \subset B_R\setminus B_r$, the existence of $\varphi_{\tilde E_n}$ follows; 
moreover, the Implicit Function Theorem applies and $\varphi_{\tilde E_n}$ is at least of class $C^{1}$. A direct computation gives
\begin{equation}\label{eqn:TangGradvarphisharpIndef}
\nabla \varphi_{\tilde{E}_{n}}  = - \frac{\left( r_0 + \varphi_{\tilde{E}_{n}} ({\theta}) \right)}{\partial_{\rho} \tilde{u}_{n}(\mathbf{x})} \nabla_T \tilde u_n
= - \frac{\left( r_0 + \varphi_{\tilde{E}_{n}} ({\theta}) \right)}{\partial_{\rho} \tilde{u}_{n}(\mathbf{x})} \left(\nabla \tilde u_n -(\nabla \tilde u_n \cdot \theta)\theta\right),
\end{equation}
where $\nabla_T$ denotes the tangential component of the gradient. Since $\tilde u_n$ converges in 
$C^{1,\alpha}$ to the radial function $w$, we deduce that $\varphi_{\tilde E_n}$ 
is of class $C^{1,\alpha}$ and that \eqref{eqn:NormTo0FavRegC1aRegSharpInde} holds true, concluding the proof.  
\end{proof}
\begin{remark}\label{rem:change_center_2}
Denote with $\bari(\tilde E_n)$ the barycenter of $\tilde E_n$. 
Thanks to Lemma \ref{lemma:favRegionAsymptSphericalIndef} we obtain that 
\[
\bari( \tilde E_n ) \to \mathbf{0}\qquad\text{as }n\to+\infty.
\]
As a consequence, by Remark 
\ref{rem:change_center_1} we have that all the results in this section, and in particular the last proposition, hold true if we center the blow-up procedure at the barycenters of the favorable regions $ E_\eps$ (and not at the unique maximum points of $ u_\eps$).
\end{remark}

To conclude this section, we use the maximum principle to show the exponential decay of $\tilde u_n$, uniformly in $n$. This will be crucial for the asymptotic location of the favorable zone, discussed in the next section.
\begin{lemma}\label{lem:exp_decay_utilde}
For every $\eta>0$,  if $n$ is sufficiently large then
\[
0\le \tilde u_n(\bx) \leq (1+\eta)w(\bx) \qquad\text{for }\bx \in \R^N \setminus (1+\eta)B.
\]
\end{lemma}
\begin{proof}
Since $\tilde u_n\equiv0$ outside $\tilde\Omega_n$, the proof is an easy application of the maximum principle in $A = \tilde\Omega_n \setminus  (1+\eta)\overline{B}$: indeed, for every $\eta>0$ and $n$ sufficiently large we have 
that 
\[
-\Delta (\tilde u_n - (1+\eta)w) + \tilde \lambda_0 \underline{m} (\tilde u_n - (1+\eta)w) = 
-(\tilde\lambda_n-\tilde\lambda_0) \underline{m} \tilde u_n \le 0\qquad\text{in }A
\]
by Lemma \ref{lemma:favRegionAsymptSphericalIndef}, and $\tilde u_n - (1+\eta)w\le 0$ on 
$\partial A$ by uniform convergence on compact sets.
\end{proof}

\section{Asymptotic location of the favorable zone}

In this section we address the asymptotic position of the favorable zone, 
i.e. the limit location of the global maximum point $\bx_n=\bx_{\eps_n}$ of the optimal eigenfunction $u_n=u_{\eps_n}$. As a byproduct of our analysis we will provide also 
an expansion of the associated eigenvalue. To proceed, we introduce the following notation:
\[
d^*\coloneqq \max_{\bp\in\Omega}\dist(\bp,\partial\Omega),
\qquad
d_n\coloneqq \dist(\bx_n,\partial\Omega) ,
\qquad
d_\infty\coloneqq \lim_{n\to+\infty} d_n
\]
(up to subsequences). Notice that, in the blow-up scale, all the above distances  are divided by $k^{1/N}$. In particular, $d_n/k^{1/N}\to+\infty$, as $n\to+\infty$, by Lemma \ref{prop:PropLimitToProblemOverRnIndef}. We will prove the following.
\begin{proposition}\label{prop:location}
Under the previous notation, we have:
\begin{enumerate}
\item $d_\infty = d^*$;
\item $\tilde\lambda_n = \tilde \lambda_0 +  e^{-2\sqrt{\tilde\lambda_0\underline{m}}(1+o(1))d^*/ k_n}$, as $n\to+\infty$.
\end{enumerate}
\end{proposition}
\begin{corollary}\label{coro:3e4di1}
Properties \ref{i:sharp_pos} and \ref{i:eig} of Theorem \ref{thm:intro_qual} hold true.
\end{corollary}

To prove the proposition, we follow the 
strategy proposed in \cite{MR1736974} to deal with the analogous problem, in the framework of 
least energy solutions to semilinear elliptic equations. As in \cite{MR1736974}, we base our analysis 
on the sharp asymptotics for the problem settled in radially symmetric domains, i.e. in a ball, 
instead of $\Omega$. We recall that in such case, by Schwarz symmetrization, the optimal favorable 
domain is a concentric ball of the appropriate measure (as proved in 
\cite[Remark 3.10]{CantrellCosner:indefiniteWeight}).
\begin{lemma}\label{lem:asballCC}
Let $B_R$ denote the ball of radius $R$, $R$ large, and $B\subset B_R$ denote the concentric ball 
of unitary volume. Then
\[
\Lambda (R) = \lambda^1(B, B_R) = \inf_{A \subset B_R : \mathcal{L}(A) = 1} \lambda^1(A, B_R)
\]
satisfies
\begin{equation}\label{eq:asympt_radial}
\Lambda(R) = \tilde \lambda_0 +  e^{-2\sqrt{\tilde\lambda_0\underline{m}}(1+o(1))R}
\qquad\text{as }r\to+\infty.
\end{equation}
\end{lemma}
\begin{proof}
This result can be proved in several different ways: an overcomplicated (although direct) 
strategy consists in deducing it as a byproduct of 
\cite[Theorem 1.3]{FerreriVerzini:AsymptSphericalIndefinite}; alternatively, one can use 
sub/supersolutions, in the spirit of Lemma \ref{lem:exp_decay_utilde} above, as done in 
\cite{MR1736974} for the semilinear case; finally, since our problem is linear, one can 
evaluate $\Lambda(R)$ as an implicit function defined in terms of Bessel functions.  
\end{proof}
Proposition \ref{prop:location} will follow by estimating 
$\tilde\lambda_n$ both from above and below, in terms of $\Lambda(R)$ evaluated at values of 
$R$ suitably depending on $d^*,d_\infty$. 
\begin{lemma}\label{lem:labo}
We have
\[
\tilde\lambda_n \le \Lambda \left(\frac{d^*}{k_n}\right).
\]
\end{lemma}
\begin{proof}
By the definitions of $d^*$  and $\tilde\Omega_n$, we have that there exists $\bx^*_n$ such that 
$B_{d^*/k_n} (\bx^*_n) \subset \tilde \Omega_n$. The lemma follows by the characterization of $\Lambda$.
\end{proof}
\begin{lemma}\label{lem:lbel}
For every $\delta>0$ sufficiently small,
\[
\tilde\lambda_n \ge \Lambda \left(\frac{d_\infty+\delta}{k_n}\right)+o\left(e^{-2\sqrt{\tilde\lambda_0\underline{m}}(d_\infty+\delta)/k_n}\right)
\qquad\text{as }n\to+\infty.
\]
\end{lemma}
\begin{proof}
Let $\delta>0$ be fixed, and $n$ sufficiently large, in such a way that  $d_\infty+\delta > d_n$. By 
definition, this means that $B_{(d_\infty+\delta)/k_n}(\mathbf{0})$ is not completely contained in 
$\tilde\Omega_n$ (and the same is true for any larger radius). On the other hand, since 
$B_{d_n/k_n}(\mathbf{0})\subsetneq\tilde\Omega_n$ (if $\tilde\Omega_n$ is a ball there is nothing to prove, by Lemma \ref{lem:asballCC}), if $\delta$ is sufficiently small we have that 
$\Lcal(B_{(d_\infty+\delta)/k_n}(\mathbf{0})) < \Lcal(\tilde \Omega_n)$. Thanks to these facts, we can choose 
$\delta'_n>\delta$ in such a way that the set 
\[
A_n=B_{(d_\infty+\delta'_n)/k_n}(\mathbf{0})\cap \tilde \Omega_n
\]
satisfies
\[
\Lcal (A_n) = \Lcal (B_{(d_\infty+\delta)/k_n}(\mathbf{0})).
\]
Notice that, as $n\to+\infty$, up to subsequences $\bx_n\to\bx_\infty$ and 
\[
\delta'_n \to \delta'>\delta,
\qquad\text{where }\qquad
\Lcal (B_{d_\infty+\delta'}(\bx_\infty)\cap  \Omega) = \Lcal (B_{d_\infty+\delta}(\mathbf{\bx_\infty})).
\]
Since $(d_\infty+\delta'_n)/k_n\to+\infty$, as $n\to +\infty$, Lemma \ref{lemma:favRegionAsymptSphericalIndef} implies that, for $n$ sufficiently large, 
$\tilde E_n\subset\subset A_n$; by \cite[Remark 3.10]{CantrellCosner:indefiniteWeight} we infer that
\begin{equation}\label{eq:confrontoCC}
\Lambda \left(\frac{d_\infty+\delta}{k_n}\right) \le 
\inf_{D \subset A_n : \mathcal{L}(D) = 1} \lambda^1(D, A_n) \le \lambda^1(\tilde E_n, A_n).
\end{equation}

We now construct an admissible function to estimate $\lambda^1(\tilde E_n, A_n)$, by 
suitably cutting-off $\tilde u_n$ on $A_n$. More precisely, let $\eta$ be a smooth cut-off function, such that
\[
\eta_n = 
\begin{cases}
1 & B_{(d_\infty+\delta'_n)/k_n-1}(\mathbf{0})\smallskip\\
0 & \R^N\setminus B_{(d_\infty+\delta'_n)/k_n}(\mathbf{0})
\end{cases}
\]
(in particular, we can assume that $|\nabla \eta_n|\le2$); 
again, by Lemma \ref{lemma:favRegionAsymptSphericalIndef}, 
$\tilde E_n\subset\subset B_{(d_\infty+\delta'_n)/k_n-1}(\mathbf{0})$ for $n$ large and 
\begin{equation}\label{eq:segnogiusto}
\eta_n(\bx)\neq1
\qquad\implies\qquad 
\tilde m_n(\bx) = -\underline{m}<0.
\end{equation}
Let us consider the function 
\[
v_n = \eta_n \tilde u_n \in H^1_0(A_n).
\]
Using \eqref{eq:segnogiusto} we obtain
\[
\int_{A_n} \tilde m_n v_n^2 = \int_{\tilde \Omega_n} \tilde m_n \eta_n^2\tilde u_n^2 =
\int_{\tilde \Omega_n} \tilde m_n \tilde u_n^2 - \int_{\tilde \Omega_n} \tilde m_n (1-\eta_n^2)\tilde u_n^2 \ge \int_{\tilde \Omega_n} \tilde m_n \tilde u_n^2 >0.
\]
Then $v$ is admissible for $\lambda^1(\tilde E_n, A_n)$ and \eqref{eq:confrontoCC} yields
\[
\begin{split}
\Lambda \left(\frac{d_\infty+\delta}{k_n}\right) &\le \frac{\int_{A_n} |\nabla v_n|^2}{\int_{A_n} \tilde m_n v_n^2}  = \frac{\int_{A_n} |\nabla (\eta_n\tilde u_n)|^2}{\int_{A_n} \tilde m_n \eta_n^2 \tilde u_n^2} = \frac{\int_{A_n} -\Delta\tilde u_n\cdot \eta_n^2 \tilde u_n + 
\int_{A_n} \tilde u_n^2|\nabla \eta_n |^2}{\int_{A_n} \tilde m_n \eta_n^2 \tilde u_n^2}\\
&= \tilde \lambda_n + \frac{\int_{A_n} \tilde u_n^2|\nabla \eta_n |^2}{\int_{A_n} \tilde m_n \eta_n^2 \tilde u_n^2}.
\end{split}
\]
Now, by strong $H^1$ convergence of $\tilde u_n$ and weak-$*$ $L^\infty$ convergence of $\tilde m_n$, we have that 
\[
\int_{A_n} \tilde m_n \eta_n^2 \tilde u_n^2 \ge \int_{\R^N} \tilde m_n  \tilde u_n^2 \to 
\int_{\R^N} \tilde m_0 w^2>0.
\]
On the other hand, recalling Lemma \ref{lem:exp_decay_utilde} and the decay of $w$ (see equation 
\eqref{eq:decay_w}), we have that, for $|\bx|$ large, 
$\tilde u_n^2(\bx) \le C r^{-(N-1)} e^{-2\sqrt{\tilde\lambda_0\underline{m}}|\bx|}$. 
Integrating, we conclude that, for $n$ sufficiently large, 
\[
\Lambda \left(\frac{d_\infty+\delta}{k_n}\right) \le \tilde \lambda_n + C e^{-2\sqrt{\tilde\lambda_0\underline{m}}(d_\infty+\delta'_n)/k_n}= 
\tilde \lambda_n + Ce^{-2\sqrt{\tilde\lambda_0\underline{m}}(1+o(1))
(d_\infty+\delta')/k_n},
\] 
and the lemma follows since $\delta'>\delta$.
\end{proof}
\begin{proof}[Proof of Proposition \ref{prop:location}]
Taking into account Lemmas \ref{lem:asballCC}, \ref{lem:lbel} and \ref{lem:labo} we can write, for every $\delta>0$ sufficiently small,
\[
\tilde \lambda_0 +  e^{-2\sqrt{\tilde\lambda_0\underline{m}}(1+o(1))(d_\infty+\delta)/k_n}\le \tilde\lambda_n \le
\tilde \lambda_0 +  e^{-2\sqrt{\tilde\lambda_0\underline{m}}(1+o(1))d^*/k_n},
\]
as $n\to+\infty$. In particular, we have $d^*\le d_\infty + \delta$, for every $\delta>0$. 
Since by definition $d_\infty\le d^*$, we deduce that actually $d_\infty= d^*$. 
Substituting in the inequality above, since $\delta$ is arbitrary we obtain the asymptotics 
of $\lambda_n$.
\end{proof}

\section{\texorpdfstring{$C^{1,1}$}{C1,1} regularity}\label{section:sharpRegularity}

In this section we provide some sharper results on the regularity of the 
eigenfunctions and of the boundary of the favourable regions. Such results are of 
interest by themselves, but they will also guarantee enough regularity to justify 
all the steps in the differentiation of the eigenvalue with respect to the shape, 
carried out in Section \ref{subsect:DifferentiabilitySharpPos}. 

We begin with a sharper analysis of the regularity of the blow-up family of eigenfunctions. To this aim, we exploit some recent results on regularity theory for transmission problems with $C^{1, \alpha}$ interfaces, proved by Caffarelli, Soria-Carro, and Stinga in \cite{Caffarelli2021:TransmissionProblems} for the case of harmonic functions, and generalized by Dong in \cite{Dong2021:TransmissionProblems} to the case of elliptic equations in divergence form.
\begin{proposition}\label{prop:C2aRegularSharpIndef}
For any $R$ sufficiently large there exists $\varepsilon_R>0$ such that, if $\varepsilon \le \varepsilon_R$, the functions $\tilde{u}_{\varepsilon}$ belong to $C^{2, \alpha}\left( \overline{\tilde{E}_{\varepsilon}} \right) \bigcap C^{2, \alpha}\left( \overline{B_R \setminus \tilde{E}_{\varepsilon}} \right)$.
\end{proposition}
\begin{proof}
The idea is to derive a transmission problem for each component of $\nabla \tilde{u}_{\varepsilon}$, and then apply the regularity result \cite[Theorem 1.2]{Dong2021:TransmissionProblems}.

Consider any partial derivative $\partial_i \tilde{u}_{\varepsilon}$ along the vectors of the orthonormal basis $\mathbf{e}_1, \dots, \mathbf{e}_N$. Applying the Gauss-Green formula and exploiting the $C^{1, \alpha}\left( \overline{ \tilde{\Omega}_{\varepsilon} } \right) \cap H^2\left(  \tilde{\Omega}_{\varepsilon} \right)$ regularity of $\tilde{u}_{\varepsilon}$ that holds by standard elliptic regularity, for any $\varphi \in C^{\infty}_c(\tilde{\Omega}_{\varepsilon})$ we can write
\begin{align}\label{eqn:eqn1C2aPropSharpIndef}
\begin{split}
    & \lambda^1(\tilde{E}_{\varepsilon}, \tilde{\Omega}_{\varepsilon}) \int_{\tilde{\Omega}_{\varepsilon}} \tilde{m}_{\varepsilon} \varphi \partial_i \tilde{u}_{\varepsilon} = 
    \lambda^1(\tilde{E}_{\varepsilon}, \tilde{\Omega}_{\varepsilon}) \left[ \overline{m} \int_{\tilde{E}_{\varepsilon}} \varphi \partial_i \tilde{u}_{\varepsilon} - 
    \underline{m} \int_{\tilde{\Omega}_{\varepsilon} \setminus \tilde{E}_{\varepsilon}} \varphi \partial_i \tilde{u}_{\varepsilon} \right] = \\
    & =\lambda^1(\tilde{E}_{\varepsilon}, \tilde{\Omega}_{\varepsilon}) (\overline{m}+\underline{m}) \int_{\partial \tilde{E}_{\varepsilon}} \tilde{u}_{\varepsilon} \varphi \nu_i - \lambda^1(\tilde{E}_{\varepsilon}, \tilde{\Omega}_{\varepsilon}) \left[ \overline{m} \int_{\tilde{E}_{\varepsilon}} \tilde{u}_{\varepsilon} \partial_i \varphi  - 
    \underline{m} \int_{\tilde{\Omega}_{\varepsilon} \setminus \tilde{E}_{\varepsilon}} \tilde{u}_{\varepsilon} \partial_i \varphi \right] ,
    \end{split}
\end{align}
where $\nu_i$ denotes the $i$-th component of the outer unit normal to $\tilde{E}_{\varepsilon}$. Using the equation of $\tilde u_\eps$ we obtain
\begin{align}\label{eqn:eqn2C2aPropSharpIndef}
\begin{split}
    & \lambda^1(\tilde{E}_{\varepsilon}, \tilde{\Omega}_{\varepsilon}) \left[ \overline{m} \int_{\tilde{E}_{\varepsilon}} \tilde{u}_{\varepsilon} \partial_i \varphi  - 
    \underline{m} \int_{\tilde{\Omega}_{\varepsilon} \setminus \tilde{E}_{\varepsilon}} \tilde{u}_{\varepsilon} \partial_i \varphi \right] = 
    - \int_{\tilde{E}_{\varepsilon}} \Delta \tilde{u}_{\varepsilon} \partial_i \varphi -
    \int_{\tilde{\Omega}_{\varepsilon} \setminus \tilde{E}_{\varepsilon}} \Delta \tilde{u}_{\varepsilon} \partial_i \varphi = \\
    & = \int_{\tilde{\Omega}_{\varepsilon}} \nabla \tilde{u}_{\varepsilon} \partial_i \nabla \varphi = - \int_{\tilde{\Omega}_{\varepsilon}} \nabla ( \partial_i \tilde{u}_{\varepsilon} ) \nabla \varphi .
    \end{split}
\end{align}
Let us denote $v \coloneqq \partial_i \tilde{u}_{\varepsilon}$. Combining \eqref{eqn:eqn1C2aPropSharpIndef} with \eqref{eqn:eqn2C2aPropSharpIndef}, we get that $v$ solves
\begin{equation}\label{eqn:eqn3C2aPropSharpIndef}
\begin{cases}
-\Delta v = \lambda^1(\tilde{E}_{\varepsilon}, \tilde{\Omega}_{\varepsilon}) \tilde{m}_{\varepsilon} v & \text{in } \tilde{E}_{\varepsilon} \bigcup \left( \tilde{\Omega}_{\varepsilon} \setminus \overline{\tilde{E}_{\varepsilon}} \right) , \\
[v] = 0, \quad [\partial_{\nu} v] = - \lambda^1(\tilde{E}_{\varepsilon}, \tilde{\Omega}_{\varepsilon}) (\overline{m}+\underline{m}) \tilde{u}_{\varepsilon} \nu_i & \text{on } \partial\tilde{E}_{\varepsilon} , \\
v = \partial_i \tilde{u}_{\varepsilon} & \text{on }
 \partial\tilde{\Omega}_{\varepsilon} , 
\end{cases}
\end{equation}
where the symbol $[ \cdot ]$ denotes the jump operator: 
\begin{equation}\label{eq:jumpop}
[ f ](\bx) = \lim_{\substack{{\by\to\bx}\\{\by\in\tilde{E}_{\varepsilon}}}} f(\by)
- \lim_{\substack{{\by\to\bx}\\{\by\in\tilde{\Omega}_{\varepsilon} \setminus \overline{\tilde{E}_{\varepsilon}}}}} f(\by)
\end{equation}
and the boundary condition is understood in the 
sense of traces. Now, for $\eps$ sufficiently small we have that $B_R\subset\subset 
\tilde\Omega_\eps$. To get rid of the boundary conditions, we notice that, since $\partial_i 
\tilde{u}_{\varepsilon} \in C^{\infty}(N_{\delta}(B_R))$ for some $\delta$ sufficiently small, 
depending on $R$ but independent of $\varepsilon$, by smooth cutoff there exists an extension 
$\bar{v}$ of $\partial_i \tilde{u}_{\varepsilon}$ in $B_R$ of class $C^{\infty}(\overline{B_R})$. 
Hence, we  can rewrite  problem \eqref{eqn:eqn3C2aPropSharpIndef} with respect to the function 
$v - \bar{v}$ in $B_R$, and apply \cite[Theorem 1.2]{Dong2021:TransmissionProblems}:  since
$\partial \tilde{E}_{\varepsilon}$ is $C^{1,\alpha}$, by Proposition \ref{prop:FavRegC1aRegSharpIndef}, we infer that $v$ is $C^{1,\alpha}$ on both sides and up to 
$\partial \tilde E_\eps$, and the proposition follows.
\end{proof}
From Proposition \ref{prop:C2aRegularSharpIndef}, we can deduce the following result.
\begin{corollary}\label{prop:C11RegularSharpIndef}
For any $R$ sufficiently large there exists $\varepsilon_R>0$ such that, if $\varepsilon \le \varepsilon_R$, the functions $\tilde{u}_{\varepsilon}$ belong to $ C^{1, 1}\left( \overline{B_R} \right)$.
\end{corollary}
\begin{proof}
Since $\tilde u_\eps \in C^{2, \alpha}\left( \overline{\tilde{E}_{\varepsilon}} \right) \bigcap C^{2, \alpha}\left( \overline{B_R \setminus \tilde{E}_{\varepsilon}} \right)$, $\tilde u_\eps \in 
C^{1, \alpha}\left(  \overline{B_R } \right)$ and $\partial \tilde E_\eps$ is of class 
$C^{1, \alpha}$, we deduce that $\tilde u_\eps \in W^{2,\infty}(B_R)$. In particular, each 
partial derivative is in $W^{1,\infty}(B_R)$ and thus it is Lipschitz continuous (see e.g. 
\cite[Section 11.4]{MR2527916}).
\end{proof}

The  next result concerns the convergence, in a certain sense, of $\tilde{u}_{\varepsilon}$ to $w$ 
in $C^{1, 1}$, as $\varepsilon$ approaches $0^{+}$.

We remark that, in general, we cannot expect a $C^{1, 1}$ convergence in $\overline{B_R}$, since 
the second derivatives of both $\tilde{u}_{\varepsilon}$ and $w$ have jumps, and the jump points 
coincide only in the limit $\varepsilon \to 0^{+}$ but never coincide for $\varepsilon \neq 0$. Indeed, recall that $\tilde{E}_{\varepsilon}$ is never exactly a ball, unless $\Omega$ is a ball \cite{Lamboley:OptimizersRobin}, while $w$ is radially symmetric by \cite[Proposition 2.1]{FerreriVerzini:AsymptSphericalIndefinite}.

Notice however that, at least for $w$, the jump occurs only for the derivative along the normal to 
the optimal region $\partial B$, while the tangential derivatives vanish. Hence, we expect that, also for the blow-up solutions $\tilde{u}_{\varepsilon}$, the tangential second derivatives to spheres vanish uniformly, where they are defined.

To make our ideas rigorous, we exploit once again the aforementioned regularity theory for interface problems. Moreover, for $\eps$ small, we consider the nearly spherical representation $\varphi_{\tilde{E}_{\varepsilon}}$ of 
$\partial \tilde{E}_{\varepsilon}$, which is of class $C^{1,\alpha}$, $\alpha<1$, by Proposition \ref{prop:FavRegC1aRegSharpIndef}. We consider the following 
family of diffeomorphisms $\mathbf{\Phi}_{\varepsilon}:\R^N \to \R^N$, written in polar coordinates $\bx = \rho \, \theta$ as
\[
\mathbf\Phi_\eps(\rho \, \theta) \coloneqq \left[ \rho^N + h(\rho) \left(( r_0+\varphi_{\tilde{E}_{\varepsilon}}(\theta) )^{N} - r_0^N \right) \right]^{1/N} \, \theta \quad \text{ for } \theta \in \sphere^{N-1}, \, \rho\ge0,
\]
where the smooth cut-off function $0\le h(\rho)\le 1$ is such that $h \equiv 1$ for $3r_0/4 \le \rho \le 5 r_0/4$ and $h \equiv 0$ outside $(r_0/2 , 3 r_0/2)$. Then one can directly check that 
$\mathbf{\Phi}_{\varepsilon}$ satisfies the following properties: 
\begin{align}\label{eqn:DiffPropertC11ConvSharpIndef}
    \begin{split}
        & \|\mathbf{\Phi}_{\varepsilon} - \mathbf{I}  \|_{C^{1, \alpha}(\R^N)} = o(1) \; \text{for some } 0 < \alpha < 1 , \quad\\
        & \mathbf{\Phi}_{\varepsilon} = \mathbf{I} \text{ in }\R^N\setminus B_{R/2},\quad\mathbf{\Phi}_{\varepsilon}(B)= \tilde{E}_{\varepsilon} , \quad
        \mathbf{\Phi}_{\varepsilon}(\partial B)= \partial \tilde{E}_{\varepsilon} ,
    \end{split}
\end{align}
where we are assuming $R>3r_0$,  and the quantities $o(1)$ are to be intended for $\varepsilon \to 0^{+}$. For more details on the properties of 
$\mathbf{\Phi}_{\varepsilon}$ we refer to 
Section \ref{section:DeformationPathsProperties}, where a more general class of such maps is treated.

Given a family of orthonormal coordinate axes $\mathbf{e}_1, \dots, \mathbf{e}_N$, we define the transported functions
\[
v_{\varepsilon, i} \coloneqq \left( \partial_i \tilde{u}_{\varepsilon} \right) \circ \mathbf{\Phi}_{\varepsilon} , \qquad i=1, \dots, N
\]
then we can prove the following.
\begin{proposition}\label{prop:TangDerVanishSharpIndef}
For any $R$ sufficiently large, in the limit $\varepsilon \to 0^{+}$ it holds that 
\[
v_{\varepsilon, i} \to \partial_i w \qquad \text{in } C^{1, \alpha}\left( \overline{B} \right) \cap C^{1, \alpha}\left( \overline{B_R \setminus B} \right)
\]
for any $i = 1, \dots, N$ and any $0 \le \alpha < 1$.
\end{proposition}
\begin{proof}
Fix an $i \in \{1, \dots, N\}$. Since $\partial_i \tilde{u}_{\varepsilon}$ solves 
\eqref{eqn:eqn3C2aPropSharpIndef}, the function $v_{\varepsilon, i}$ solves in $H^1(B_R)$
\begin{equation}\label{eqn:eqn1TangDerVanishSharpIndef}
\begin{cases}
-\diverg(A_{\varepsilon} \nabla v_{\varepsilon, i}) = \lambda^1(\tilde{E}_{\varepsilon}, \tilde{\Omega}_{\varepsilon}) \tilde{m}_{0} v_{\varepsilon, i} J_{\varepsilon} & \text{in } B \bigcup \left( B_R \setminus \overline{B} \right) , \\
[v_{\varepsilon, i}] = 0, \quad [\partial_{\mathbf{n}} v_{\varepsilon, i}] = - \lambda^1(\tilde{E}_{\varepsilon}, \tilde{\Omega}_{\varepsilon}) (\overline{m}+\underline{m}) (\tilde{u}_{\varepsilon} \nu_i) \circ \mathbf{\Phi}_{\varepsilon} J_{\varepsilon, T} & \text{on } \partial B, \\
v_{\varepsilon, i} = \left( \partial_i \tilde{u}_{\varepsilon} \right) \circ \mathbf{\Phi}_{\varepsilon} & \text{on }
 \partial B_R , 
\end{cases}
\end{equation}
where $\mathbf{n}$ and $\nu$ are the unitary outer normal vectors to $\partial B$ and $\partial \tilde{E}_{\varepsilon}$ respectively, while
\begin{align}\label{eqn:GeometricC1aDefSharpIndef}
\begin{split}
J_{\varepsilon} \coloneqq \vert \det & \left( D \mathbf{\Phi}_{\varepsilon} \right) \vert , \quad A_{\varepsilon} \coloneqq J_{\varepsilon} D\mathbf{\Phi}_{\varepsilon}^{-1} D\mathbf{\Phi}_{\varepsilon}^{-T} \qquad \text{in } \R^N \\
& J_{\varepsilon, T} \coloneqq \vert D\mathbf{\Phi}_{\varepsilon}^{-T} \mathbf{n} \vert J_{\varepsilon} \qquad \text{on } \partial B .
\end{split}
\end{align}
Taking into account the fact that the function $v_{0, i} \coloneqq \partial_i w$ solves in $H^1(B_R)$
\begin{equation}\label{eqn:eqn2TangDerVanishSharpIndef}
\begin{cases}
-\Delta v_{0, i} = \lambda^1(B, \R^N) \tilde{m}_{0} v_{0, i} & \text{in } B \bigcup \left( B_R \setminus \overline{B} \right) , \\
[v_{0, i}] = 0, \quad [\partial_{\mathbf{n}} v_{0, i}] = - \lambda^1(B, \R^N) (\overline{m}+\underline{m}) w n_i  & \text{on } \partial\tilde{E}_{\varepsilon} , \\
v_{0, i} = \partial_i w  & \text{on }
 \partial B_R , 
\end{cases}
\end{equation}
combining \eqref{eqn:eqn1TangDerVanishSharpIndef} and \eqref{eqn:eqn2TangDerVanishSharpIndef} we can write the following problem, in $H^1(B_R)$, for the difference $z_{\varepsilon, i} \coloneqq v_{\varepsilon, i} - v_{0, i}$:
\begin{equation}\label{eqn:eqn3TangDerVanishSharpIndef}
\begin{cases}
-\diverg(A_{\varepsilon} \nabla z_{\varepsilon, i}) = \diverg(F_{\varepsilon, i}) + f_{\varepsilon, i} & \text{in } B \bigcup \left( B_R \setminus \overline{B} \right) , \\
[z_{\varepsilon, i}] = 0, \quad [\partial_{\mathbf{n}} z_{\varepsilon, i}] = g_{\varepsilon, i} & \text{on } \partial B, \\
z_{\varepsilon, i} = h_{\varepsilon, i} & \text{on }
 \partial B_R , 
\end{cases}
\end{equation}
where we have defined
\begin{align}\label{eqn:DefC1aVanishQuantSharpIndef}
\begin{split}
    & F_{\varepsilon, i} \coloneqq (A_{\varepsilon} - \mathbf{I}) \nabla v_{0, i} , \qquad f_{\varepsilon, i} \coloneqq \lambda^1(\tilde{E}_{\varepsilon}, \tilde{\Omega}_{\varepsilon}) \tilde{m}_{0} v_{\varepsilon, i} J_{\varepsilon} - \lambda^1(B, \R^N) \tilde{m}_{0} v_{0, i} , \\
    & g_{\varepsilon, i} \coloneqq (\overline{m}+\underline{m}) \left[ \lambda^1(B, \R^N) w n_i - \lambda^1(\tilde{E}_{\varepsilon}, \tilde{\Omega}_{\varepsilon}) (\tilde{u}_{\varepsilon} \nu_i) \circ \mathbf{\Phi}_{\varepsilon} J_{\varepsilon, T} \right] , \\
    & h_{\varepsilon, i} \coloneqq \partial_i \tilde{u}_{\varepsilon} - \partial_i w \text{ on }\partial B_R
\end{split}
\end{align}
(recall that $\mathbf{\Phi}$ is the identity near $\partial B_R$). 
Now we wish to reduce problem \eqref{eqn:eqn3TangDerVanishSharpIndef} to the case of homogeneous 
Dirichlet boundary conditions. To this aim, we recall that, by standard (interior) elliptic 
regularity, it holds that $h_{\varepsilon, i} = \partial_i \tilde{u}_{\varepsilon} 
-  \partial_i w \to 0$, in a neighborhood of $\partial B_R$, in $C^k$ for every $k\ge0$. 
Using a smooth  
cutoff, we can easily construct functions  $\bar{z}_{\varepsilon, i} \in C^{1, \alpha}\left( \overline{B_R} \right)$ satisfying
\begin{equation}\label{eqn:eqn9TangDerVanishSharpIndef}
\bar{z}_{\varepsilon, i} =  h_{\varepsilon, i}\text{ on }\partial B_R, \quad 
\bar{z}_{\varepsilon, i} =  0\text{ in } B_{(1-\eta)R}, \quad \| \bar{z}_{\varepsilon, i} \|_ {C^{1, \alpha}\left( \overline{B_R} \right)} \to 0 \text{ for } \varepsilon \to 0^+.
\end{equation}
Exploiting $\bar{z}_{\varepsilon, i}$ and introducing the quantity
\begin{equation}\label{eqn:2DefC1aVanishQuantSharpIndef}
    \tilde{F}_{\varepsilon, i} \coloneqq F_{\varepsilon, i} + A_{\varepsilon} \nabla \bar{z}_{\varepsilon, i} ,
\end{equation}
we can rewrite \eqref{eqn:eqn3TangDerVanishSharpIndef} as
\begin{equation}\label{eqn:eqn10TangDerVanishSharpIndef}
\begin{cases}
-\diverg(A_{\varepsilon} \nabla (z_{\varepsilon, i} - \bar{z}_{\varepsilon, i})) = \diverg(\tilde{F}_{\varepsilon, i}) + f_{\varepsilon, i} & \text{in } B \bigcup \left( B_R \setminus \overline{B} \right) , \\
[z_{\varepsilon, i}] = 0, \quad [\partial_{\mathbf{n}} z_{\varepsilon, i}] = g_{\varepsilon, i} & \text{on } \partial B, \\
z_{\varepsilon, i} - \bar{z}_{\varepsilon, i} = 0 & \text{on } \partial B_R , 
\end{cases}
\end{equation}
Now we are in the position to use again \cite[Theorem 1.2]{Dong2021:TransmissionProblems} to deduce that
\begin{align}\label{eqn:eqn11TangDerVanishSharpIndef}
\begin{split}
& \| z_{\varepsilon, i} - \bar{z}_{\varepsilon, i} \|_{C^{1, \alpha} \left( \overline{B} \right)} + \| z_{\varepsilon, i} - \bar{z}_{\varepsilon, i} \|_{C^{1, \alpha} \left( \overline{B_R \setminus B} \right)} \le \\
& \le C\left( \| g_{\varepsilon, i} \|_{C^{0, \alpha}(\partial B)} +  \| \tilde{F}_{\varepsilon, i} \|_{C^{0, \alpha}(B)} +  \| \tilde{F}_{\varepsilon, i} \|_{C^{0, \alpha}\left( \overline{B_R \setminus B} \right)} +  \| f_{\varepsilon, i} \|_{L^{\infty}(B_R)} \right).
\end{split}
\end{align}
for some positive constant $C$ independent of $\varepsilon$, thanks to \eqref{eqn:DiffPropertC11ConvSharpIndef}.
Thus, thanks to \eqref{eqn:eqn9TangDerVanishSharpIndef}, to conclude the proof of Proposition 
\ref{prop:TangDerVanishSharpIndef} it is sufficient to establish that the right hand side of \eqref{eqn:eqn11TangDerVanishSharpIndef} is infinitesimal for $\varepsilon \to 0^+$. 

For the reader's convenience, we split this estimate in Lemmas \ref{lemma:TildeFC1aTo0SharpIndef}, \ref{lemma:fC1aTo0SharpIndef}, \ref{lemma:gC1aTo0SharpIndef} below, which are preceded by the preliminary Lemma \ref{lemma:GeometricC1aConvSharpIndef} dealing with the properties of geometrical quantities.
\end{proof}
\begin{lemma}\label{lemma:GeometricC1aConvSharpIndef}
Given the definitions in \eqref{eqn:GeometricC1aDefSharpIndef}, it holds that
\begin{itemize}
    \item[(i)] $\| J_{\varepsilon} - 1 \|_{C^{0, \alpha}(\R^N)} \to 0$ for $\varepsilon \to 0^+$,
    \item[(ii)] $\| A_{\varepsilon} - \mathbf{I} \|_{C^{0, \alpha}(\R^N)} \to 0$ for $\varepsilon \to 0^+$,    
    \item[(iii)] $\| J_{\varepsilon, T} - 1 \|_{C^{0, \alpha}(\partial B)} \to 0$ for $\varepsilon \to 0^+$,
    \item[(iv)] $\| \nu \circ \mathbf{\Phi}_{\varepsilon} - \mathbf{n} \|_{C^{0, \alpha}(\partial B)} \to 0$ for $\varepsilon \to 0^+$,
\end{itemize}
where $\nu$ and $\mathbf{n}$ are the unitary outer normal vectors to $\partial \tilde{E}_{\varepsilon}$ and $\partial B$, respectively.
\end{lemma}
\begin{proof}
    Let us prove point (i). The fact that $\| J_{\varepsilon} - 1 \|_{C^{0}(\R^N)} \to 0$ for $\varepsilon \to 0^+$ follows directly from \eqref{eqn:DiffPropertC11ConvSharpIndef}. Moreover, since the map $\det : M^{N \times N} \to \R$ is of class $C^{\infty}$, by composition and exploiting again \eqref{eqn:DiffPropertC11ConvSharpIndef}, it follows that $ \| \vert \det \left( D \mathbf{\Phi}_{\varepsilon} \right) \vert - 1 \|_{C^{0, \alpha}(\R^N)} = \| \det \left( D \mathbf{\Phi}_{\varepsilon} \right) - 1 \|_{C^{0, \alpha}(\R^N)} \to 0$ for $\varepsilon \to 0^+$.
    
    For point (ii) we notice that, thanks to \eqref{eqn:DiffPropertC11ConvSharpIndef}, it holds that $\| D\mathbf{\Phi}_{\varepsilon}^{-1} - \mathbf{I} \|_{C^{0, \alpha}(\R^N)} \to 0$ for $\varepsilon \to 0^+$. An analogous property holds for $D\mathbf{\Phi}_{\varepsilon}^{-T}$. Hence, exploiting point $(i)$ the conclusion follows.
    
    To deal with point (iii), we first notice that, thanks to \eqref{eqn:DiffPropertC11ConvSharpIndef}, we have $\| D\mathbf{\Phi}_{\varepsilon}^{-T} \mathbf{n} - \mathbf{n} \|_{C^{0, \alpha}(\partial B)} \to 0$ for $\varepsilon \to 0^+$. In a similar fashion, it also holds that $\| \vert D\mathbf{\Phi}_{\varepsilon}^{-T} \mathbf{n} \vert - 1 \|_{C^{0, \alpha}(\partial B)} = \| \vert D\mathbf{\Phi}_{\varepsilon}^{-T} \mathbf{n} \vert - \vert \mathbf{n} \vert \|_{C^{0, \alpha}(\partial B)} \to 0$ for $\varepsilon \to 0^+$. The conclusion follows exploiting point $(i)$.

    Now we study point (iv). We recall that
    \[
    \nu \circ \mathbf{\Phi}_{\varepsilon} = \frac{D\mathbf{\Phi}_{\varepsilon}^{-T} \mathbf{n}}{\vert D\mathbf{\Phi}_{\varepsilon}^{-T} \mathbf{n} \vert} ,
    \]
    so that the proof immediately follows from the properties introduced in the analysis of point (iii).
\end{proof}
\begin{lemma}\label{lemma:TildeFC1aTo0SharpIndef}
Let $\tilde{F}_{\varepsilon, i}$ be defined as in \eqref{eqn:2DefC1aVanishQuantSharpIndef}. Then
\[
\| \tilde{F}_{\varepsilon, i} \|_{C^{0, \alpha}(B_R)} \to 0 \qquad \text{for } \varepsilon \to 0^+ .
\]
\end{lemma}
\begin{proof}
Since by definition $\tilde{F}_{\varepsilon, i} = F_{\varepsilon, i} + A_{\varepsilon} \nabla \bar{z}_{\varepsilon, i}$, where $F_{\varepsilon, i}$ is defined in \eqref{eqn:DefC1aVanishQuantSharpIndef}, thanks to \eqref{eqn:eqn9TangDerVanishSharpIndef} and Lemma \ref{lemma:GeometricC1aConvSharpIndef} (ii) it is sufficient to prove that $\| F_{\varepsilon, i} \|_{C^{0, \alpha}(B)} \to 0$ for $\varepsilon \to 0^+$. This follows immediately from Lemma \ref{lemma:GeometricC1aConvSharpIndef} (ii).
\end{proof}
\begin{lemma}\label{lemma:fC1aTo0SharpIndef}
Let $f_{\varepsilon, i}$ be defined as in \eqref{eqn:DefC1aVanishQuantSharpIndef}. Then
\[
\| f_{\varepsilon, i} \|_{L^{\infty}(B_R)} \to 0 \qquad \text{for } \varepsilon \to 0^+ .
\]
\end{lemma}
\begin{proof}
Let us decompose $f_{\varepsilon, i} = f_{1} + f_2 + f_3$, where
\begin{align*}
f_1 \coloneqq \left[ \lambda^1(\tilde{E}_{\varepsilon}, \tilde{\Omega}_{\varepsilon})  -  \lambda^1( B, \R^N) \right] & \tilde{m}_{0} v_{\varepsilon, i} J_{\varepsilon} , \qquad
f_2 \coloneqq \lambda^1(B, \R^N) \tilde{m}_{0} v_{\varepsilon, i} \left[ J_{\varepsilon} -1 \right] , \\
f_3 \coloneqq & \lambda^1(B, \R^N) \tilde{m}_{0} \left[ v_{\varepsilon, i} - v_{0, i} \right] .
\end{align*}
Now,  for $\varepsilon \to 0^+$, $\| f_1 \|_{L^{\infty}(B_R)} \to 0 $ by Lemma \ref{lem:conv2la_0Indef}, while  $\| f_2 \|_{L^{\infty}(B_R)} \to 0 $ by Lemma \ref{lemma:GeometricC1aConvSharpIndef} (i).
For $f_3$ we proceed as follows. We write
\[
v_{\varepsilon, i} - v_{0, i} = (\partial_i \tilde{u}_{\varepsilon} ) \circ \mathbf{\Phi}_{\varepsilon} - \partial_i w = \left( \partial_i \tilde{u}_{\varepsilon} - \partial_i w  \right) \circ \mathbf{\Phi}_{\varepsilon} + \left((\partial_i w ) \circ \mathbf{\Phi}_{\varepsilon} - \partial_i w \right) ,
\]
and notice that $\| \left( \partial_i \tilde{u}_{\varepsilon} - \partial_i w  \right) \circ \mathbf{\Phi}_{\varepsilon} \|_{L^{\infty}(B_R)} \to 0 $ for $\varepsilon \to 0^+$ thanks to Lemma  \ref{lem:H1strong}, while $\| (\partial_i w ) \circ \mathbf{\Phi}_{\varepsilon} - \partial_i w \|_{L^{\infty}(B_R)} \to 0 $ for $\varepsilon \to 0^+$ thanks to the $C^{1, \alpha}\left( \overline{B_R} \right)$ regularity of $w$ and \eqref{eqn:DiffPropertC11ConvSharpIndef}.
\end{proof}
\begin{lemma}\label{lemma:gC1aTo0SharpIndef}
Let $g_{\varepsilon, i}$ be defined as in \eqref{eqn:DefC1aVanishQuantSharpIndef}. Then
\[
\| g_{\varepsilon, i} \|_{C^{0, \alpha}(\partial B)} \to 0 \qquad \text{for } \varepsilon \to 0^+ .
\]
\end{lemma}
\begin{proof}
We can decompose $g_{\varepsilon, i} = g_1 + g_2 + g_3 + g_4$, where
\begin{align*}
& g_1 \coloneqq \left[ \lambda^1( B, \R^N) - \lambda^1(\tilde{E}_{\varepsilon}, \tilde{\Omega}_{\varepsilon}) \right] (\overline{m}+\underline{m}) (\tilde{u}_{\varepsilon} \nu_i) \circ \mathbf{\Phi}_{\varepsilon} J_{\varepsilon, T} , \\
& g_2 \coloneqq \lambda^1(B, \R^N) (\overline{m}+\underline{m}) (\tilde{u}_{\varepsilon} \nu_i) \circ \mathbf{\Phi}_{\varepsilon} \left[ 1 -J_{\varepsilon, T} \right] , \\
& g_3 \coloneqq \lambda^1(B, \R^N) (\overline{m}+\underline{m}) \tilde{u}_{\varepsilon} \circ \mathbf{\Phi}_{\varepsilon} \left[ n_i - \nu_i \circ \mathbf{\Phi}_{\varepsilon} \right] , \\
& g_3 \coloneqq \lambda^1(B, \R^N) (\overline{m}+\underline{m}) n_i \left[ w - \tilde{u}_{\varepsilon} \circ \mathbf{\Phi}_{\varepsilon} \right] .
\end{align*}
The facts that $\|g_1\|_{C^{0, \alpha}}(\partial B) \to 0$, $\|g_2\|_{C^{0, \alpha}}(\partial B) \to 0$ and $\|g_3\|_{C^{0, \alpha}}(\partial B) \to 0$ for $\varepsilon \to 0^+$ follow readily from Lemma \ref{lem:conv2la_0Indef}, Lemma \ref{lemma:GeometricC1aConvSharpIndef} (iii) and (iv), respectively. For what concerns $g_4$, we notice that 
\[
\|\tilde{u}_{\varepsilon} \circ \mathbf{\Phi}_{\varepsilon} - w \|_{C^{0, \alpha}\left( \overline{B_R} \right)} \le \|\tilde{u}_{\varepsilon} \circ \mathbf{\Phi}_{\varepsilon} - w \circ \mathbf{\Phi}_{\varepsilon} \|_{C^{0, \alpha}\left( \overline{B_R} \right)} + \|w \circ \mathbf{\Phi}_{\varepsilon} - w \|_{C^{0, \alpha}\left( \overline{B_R} \right)}
\]
and the right hand side vanishes for $\varepsilon \to +\infty$, thanks to Lemma  \ref{lem:H1strong} and \eqref{eqn:DiffPropertC11ConvSharpIndef}.
\end{proof}

Once Proposition \ref{prop:TangDerVanishSharpIndef} is established, we 
draw the following direct consequence.
\begin{corollary}\label{corollary:UeSmallTangDerivSharpIndef}
For any $r<R$, with $B_r\subset\subset B \subset\subset B_R$, in the limit $\varepsilon \to 0^+$ it holds that
\[
   \| \nabla \tilde{u}_{\varepsilon} - (\nabla \tilde{u}_{\varepsilon} \cdot \mathbf{n} ) \mathbf{n}  \|_{C^{1, \alpha}\left( \overline{\tilde{E}_{\varepsilon} \setminus B_r} \right)} + 
    \| \nabla \tilde{u}_{\varepsilon} - (\nabla \tilde{u}_{\varepsilon} \cdot \mathbf{n} ) \mathbf{n} \|_{C^{1, \alpha}\left( \overline{B_R \setminus \tilde{E_{\varepsilon}}} \right)} \to 0 ,
\]
where $\mathbf{n}=\dfrac{\bx}{|\bx|}$ denotes the unitary outer normal vector to $B$, extended by radial prolongation over $\R^N\setminus\{\mathbf{0}\}$. 
\end{corollary}
\begin{proof}
First, we claim that
\[
\| \left(\nabla \tilde{u}_{\varepsilon} - (\nabla \tilde{u}_{\varepsilon} \cdot \mathbf{n} ) \mathbf{n} \right) \circ \mathbf{\Phi}_{\varepsilon} \|_{C^{1, \alpha}\left( \overline{B \setminus B_r} \right)} + 
    \| \left(\nabla \tilde{u}_{\varepsilon} - (\nabla \tilde{u}_{\varepsilon} \cdot \mathbf{n} ) \mathbf{n} \right) \circ \mathbf{\Phi}_{\varepsilon} \|_{C^{1, \alpha}\left( \overline{B_R \setminus B} \right)} \to 0.
\]
Indeed, it is sufficient to notice that, by Proposition 
\ref{prop:TangDerVanishSharpIndef} it holds that
\begin{align*}
& \| \left(\nabla \tilde{u}_{\varepsilon} - (\nabla  \tilde{u}_{\varepsilon} \cdot \mathbf{n} )\mathbf{n} \right) \circ \mathbf{\Phi}_{\varepsilon} - \left( \nabla w - (\nabla w \cdot \mathbf{n} ) \mathbf{n} \right)\|_{C^{1, \alpha}\left( \overline{B \setminus B_r} \right)} + \\
& +
\| \left(\nabla \tilde{u}_{\varepsilon} - (\nabla \tilde{u}_{\varepsilon} \cdot \mathbf{n} ) \mathbf{n} \right) \circ \mathbf{\Phi}_{\varepsilon} - \left( \nabla w - (\nabla w \cdot \mathbf{n} ) \mathbf{n} \right) \|_{C^{1, \alpha}\left( \overline{B_R \setminus B} \right)} \to 0
\end{align*}
for $\varepsilon \to 0^+$, and that, since $w$ is radially symmetric, it also holds that
\[
\nabla w - (\nabla w \cdot \mathbf{n} ) \mathbf{n} \equiv \mathbf{0} \quad \text{in } \overline{B_R \setminus B_r} .
\]

Once the claim is proved, to conclude just notice that 
\[
 \nabla \tilde{u}_{\varepsilon} - (\nabla \tilde{u}_{\varepsilon} \cdot \mathbf{n} ) \mathbf{n} = \left[ \left(\nabla \tilde{u}_{\varepsilon} - (\nabla  \tilde{u}_{\varepsilon} \cdot \mathbf{n}) \mathbf{n} \right) \circ \mathbf{\Phi}_{\varepsilon} \right] \circ \mathbf{\Phi}_{\varepsilon}^{-1} ,
\]
that $\| \mathbf{\Phi}_{\varepsilon}^{-1} - \mathbf{I} \|_{C^{1, \alpha}} \to 0$ for $\varepsilon \to 0^+$ and that $\mathbf{\Phi}_{\varepsilon}(\partial B) = \partial \tilde{E}_{\varepsilon}$ (see \eqref{eqn:DiffPropertC11ConvSharpIndef}).
\end{proof}

Once the $C^{1,1}$ regularity of the eigenfunctions is established, we turn to 
the regularity of the free boundary, i.e.\ of its nearly spherical representation 
$\varphi_{\tilde E_\eps}$. 

Given any $\theta_1$, $\theta_2 \in  \sphere^{N-1}$ let us give the following definitions:
\[
\mathbf{x}_i \coloneqq (r_0 + \varphi_{\tilde{E}_{\varepsilon}}(\theta_i) ) \theta_i, \qquad \text{for } i = 1, 2,
\]
$\gamma_{\theta_1, \theta_2}: [0, L(\gamma_{\theta_1, \theta_2})] \to \sphere^{N-1}$ denotes a geodesic on $\sphere^{N-1}$ connecting $\theta_1$ with $\theta_2$, with velocity $1$, and $\Gamma_{\theta_1, \theta_2}: [0, L(\gamma_{\theta_1, \theta_2})] \to \partial \tilde{E}_{\varepsilon}$ is defined as
\[
\Gamma_{\theta_1, \theta_2} \coloneqq (r_0 + \varphi_{\tilde{E}_{\varepsilon}} \circ \gamma_{\theta_1, \theta_2} ) \gamma_{\theta_1, \theta_2} .
\]

As a consequence of Proposition \ref{prop:FavRegC1aRegSharpIndef}, we can prove the following
\begin{lemma}\label{corollary:AsymptSphereGammaSharpIndef}
There exists $\bar{\varepsilon}>0$ such that, for any $0<\varepsilon \le \bar{\varepsilon}$ and any 
$\theta_1, \theta_2 \in \sphere^{N-1}$,
\begin{itemize}
    \item[(i)] $L(\Gamma_{\theta_1, \theta_2}) = (1+o(1) )L(\gamma_{\theta_1, \theta_2})$,
    \item[(ii)] $\vert \mathbf{x}_1 - \mathbf{x}_2 \vert = (1+o(1)) \vert \theta_1 - \theta_2 \vert$,
\end{itemize}
where all $o(1)$ terms are intended for $\varepsilon \to 0^{+}$ and uniformly in $\theta_1,\theta_2$.
\end{lemma}
\begin{proof}
We begin with point (i). First of all, it holds that
\[
\dot \Gamma_{\theta_1, \theta_2} = (r_0 + \varphi_{\tilde{E}_{\varepsilon}} \circ \gamma_{\theta_1, \theta_2} ) \dot \gamma_{\theta_1, \theta_2} + (\nabla_T \varphi_{\tilde{E}_{\varepsilon}}  \cdot \dot \gamma_{\theta_1, \theta_2}) \gamma_{\theta_1, \theta_2} ,
\]
so that
\[
    \vert r_0 + \varphi_{\tilde{E}_{\varepsilon}} \circ \gamma_{\theta_1, \theta_2}  \vert  - \vert \nabla_T \varphi_{\tilde{E}_{\varepsilon}}  \cdot \dot \gamma_{\theta_1, \theta_2} \vert \le \vert \dot \Gamma_{\theta_1, \theta_2} \vert 
    \le  \vert r_0 + \varphi_{\tilde{E}_{\varepsilon}} \circ \gamma_{\theta_1, \theta_2}  \vert + \vert \nabla_T \varphi_{\tilde{E}_{\varepsilon}}  \cdot \dot \gamma_{\theta_1, \theta_2} \vert
\]
and thanks to Proposition \ref{prop:FavRegC1aRegSharpIndef} we can write
$\vert \dot \Gamma_{\theta_1, \theta_2} \vert = 1 + o(1)$ . Point (i) follows immediately.

Now we turn to point (ii). We have
\[
    \vert \theta_1 - \theta_2 \vert -  \vert \varphi_{\tilde{E}_{\varepsilon}}(\theta_1) \theta_1 - \varphi_{\tilde{E}_{\varepsilon}}(\theta_2) \theta_2  \vert \le \vert \mathbf{x}_1 - \mathbf{x}_2 \vert \le \vert \theta_1 - \theta_2 \vert + \vert \varphi_{\tilde{E}_{\varepsilon}}(\theta_1) \theta_1 - \varphi_{\tilde{E}_{\varepsilon}}(\theta_2) \theta_2 \vert
\]
\and the result follows again from Proposition \ref{prop:FavRegC1aRegSharpIndef}.
\end{proof}

Now, exploiting the regularity results previously obtained 
for $\tilde{u}_{\varepsilon}$, we can extend Proposition \ref{prop:FavRegC1aRegSharpIndef} to the case $\alpha = 1$.
\begin{proposition}\label{prop:FavRegC11RegSharpIndef}
There exists $\bar{\varepsilon}>0$ such that, for any $0<\varepsilon \le \bar{\varepsilon}$, it holds that $\partial \tilde{E}_{\varepsilon}$ is nearly spherical of class $C^{1, 1}$ and
\[
\| \varphi_{\tilde{E}_{\varepsilon}} \|_{C^{1, 1}} \to 0 ,\qquad \text{for } \varepsilon \to 0^+ .
\]
\end{proposition}
\begin{proof}
Given the $C^{1, 1}$ regularity of $\tilde{u}_{\varepsilon}$ stated in Corollary \ref{prop:C11RegularSharpIndef}, exactly as in Proposition \ref{prop:FavRegC1aRegSharpIndef} we can prove that $\partial \tilde{E}_{\varepsilon}$ is nearly spherical of class $C^{1, 1}$ and $\| \varphi_{\tilde{E}_{\varepsilon}} \|_{C^{1, \alpha}} \to 0$ for $\varepsilon \to 0^+$, for any $0 \le \alpha < 1$. Hence, the only thing left to prove is that
\begin{equation}\label{eqn:ToProveFavRegC11RegSharpIndef}
\| \nabla_T \varphi_{\tilde{E}_{\varepsilon}} \|_{C^{0, 1}} \to 0 \qquad \text{for } \varepsilon \to 0^+ .
\end{equation}

To prove \eqref{eqn:ToProveFavRegC11RegSharpIndef}, the fundamental tool will be Corollary \ref{corollary:UeSmallTangDerivSharpIndef}. Fix $r<r_0<R$. To begin with, we remark that since $\partial_n \tilde{u}_{\varepsilon} = \nabla \tilde{u}_{\varepsilon} \cdot \mathbf{n}$, thanks to Proposition \ref{prop:C2aRegularSharpIndef} and Corollary \ref{corollary:UeSmallTangDerivSharpIndef}, for $\varepsilon \to 0^+$ it holds that, for 
$\alpha<1$,
 \begin{equation}\label{eqn:NablaTC1Conv0SharpIndef}
   \left\| \frac{1}{\partial_n \tilde{u}_{\varepsilon}} \left( \nabla \tilde{u}_{\varepsilon} - (\nabla \tilde{u}_{\varepsilon} \cdot \mathbf{n} ) \mathbf{n} \right)  \right\|_{C^{1, \alpha}\left( \overline{\tilde{E}_{\varepsilon} \setminus B_r} \right)} + 
    \left\| \frac{1}{\partial_n \tilde{u}_{\varepsilon}} \left( \nabla \tilde{u}_{\varepsilon} - (\nabla \tilde{u}_{\varepsilon} \cdot \mathbf{n} ) \mathbf{n} \right) \right\|_{C^{1, \alpha}\left( \overline{B_R \setminus \tilde{E_{\varepsilon}}} \right)} \to 0 ,
   \end{equation}
Let $\mathbf{x}_i \coloneqq (r_0 + \varphi_{\tilde{E}_{\varepsilon}}(\theta_i) ) \theta_i$, 
$i = 1, 2$. Then, using Lemma \ref{corollary:AsymptSphereGammaSharpIndef}, we have
\[
\frac{|\nabla_T \varphi_{\tilde{E}_{\varepsilon}}(\theta_2)-\nabla_T \varphi_{\tilde{E}_{\varepsilon}}(\theta_1)|}{|\theta_2-\theta_1|}\leq 2\frac{|\nabla_T \varphi_{\tilde{E}_{\varepsilon}}(\theta_2)-\nabla_T \varphi_{\tilde{E}_{\varepsilon}}(\theta_1)|}{|\bx_2-\bx_1|}.
\]
Recalling equation \eqref{eqn:TangGradvarphisharpIndef} and using 
\eqref{eqn:NablaTC1Conv0SharpIndef} one can easily conclude.
\end{proof}
\begin{remark}\label{rem:change_center_3}
The nearly spherical representation $\varphi_{\tilde E_\eps}$ depends of course on 
the choice of the centers of the blow-up procedure. In view of Remark 
\ref{rem:change_center_2} we have that Proposition 
\ref{prop:FavRegC11RegSharpIndef} holds true both in case such centers are taken at 
the unique maximum point of $u_\eps$, and when they are identified with the 
barycenter of the optimal favorable region $E_\eps$.  
\end{remark}
\begin{proof}[End of the proof of Theorem \ref{thm:intro_qual}]
In view of Corollaries \ref{coro:1di1} and \ref{coro:3e4di1} we are left to prove that properties 
\ref{i:varphi} and \ref{i:asimpt_spher} in the theorem hold true, regardless of the blow-up procedure 
being centered at the maximum points of $u_\eps$ or at the barycenters of $E_\eps$. In turn, taking 
into account the definition of $\tilde E_\eps$ \eqref{eq:BlowUpSets} and the blow-up scaling \eqref{eq:kbeta}, \eqref{eq:BlowUpSets}, 
such properties are direct consequences of Proposition \ref{prop:FavRegC11RegSharpIndef} and 
Remark \ref{rem:change_center_3} above, by simply writing
\[
\varphi_{\varepsilon}=\varphi_{\tilde{E}_{\varepsilon}}.\qedhere
\]
\end{proof}

\section{Quantitative asymmetry estimates for the problem in 
\texorpdfstring{$\R^N$}{RN}}\label{sec:SharpAsymFavRegSharpIndef}

This section is mainly devoted to the proof of Theorem \ref{thm:quantitStabRNSharpPos}. Before dealing 
with it, we first notice that such result readily implies both Corollary \ref{coro:quantitative_GN} and 
Theorem \ref{thm:intro_quant_domain}.
\begin{proof}[Proof of Corollary \ref{coro:quantitative_GN}]
Using Theorem \ref{thm:quantitStabRNSharpPos} and the Gagliardo-Nirenberg inequality 
we obtain, for every $0<\alpha<1$, 
\[
\begin{split}
\| \varphi_{\mathcal{A}} \|_{C^{1,\alpha}(\sphere^{N-1})}&\le C_1 
\|\varphi_{\mathcal{A}}\|_{W^{2,\infty}(\sphere^{N-1})}^{(2+N+2\alpha)/(4+N)} 
\|\varphi_{\mathcal{A}}\|_{L^{2}(\sphere^{N-1})}^{(2-2\alpha)/(4+N)} + C_2 \|\varphi_{\mathcal{A}}\|_{L^{2}(\sphere^{N-1})} \\ 
&\le C_3 
\delta^{(2+N+2\alpha)/(4+N)} 
\|\varphi_{\mathcal{A}}\|_{L^{2}(\sphere^{N-1})}^{(2-2\alpha)/(4+N)},
\end{split}
\]
where $C_1,C_2,C_3$ only depend on $N$ and $\alpha$, see e.g.~\cite[Thm. 1]{MR208360} (actually, such result is stated on bounded domains of $\R^N$, but it can be extended to $\sphere^{N-1}$ in a standard way, by reasoning as in e.g.~\cite[Sec. 2.6]{MR1688256}).
\end{proof}
\begin{proof}[Proof of Theorem  \ref{thm:intro_quant_domain}]
Using the same notation of the previous sections, let $\tilde{E}_\eps$ denote the optimal favorable set 
in the blow-up scale, centered at  $\bari(E_\eps)$. In particular, if $\varphi_\eps$ is defined as 
in Theorem \ref{thm:intro_qual} (with $\bari(E_\eps)$ instead of $\bx_\eps$), we have that 
\[
\varphi_{\varepsilon}=\varphi_{\tilde{E}_{\varepsilon}}.
\]
Since  $\tilde{\Omega}_{\varepsilon} \subset \R^N$, we infer
\[
\lambda^1(\tilde{E}_{\varepsilon}, \tilde{\Omega}_{\varepsilon}) - \lambda^1(B, \R^N) 
\ge \lambda^1(\tilde{E}_{\varepsilon}, \R^N) - \lambda^1(B, \R^N) .
\]
Now, by assumption $\bari(\tilde{E}_{\varepsilon})=\mathbf{0}$, while 
$\| \varphi_{\tilde{E}_{\varepsilon}} \|_{C^{1, 1}(\sphere^{N-1})}$ can be made arbitrarily 
small using Proposition \ref{prop:FavRegC11RegSharpIndef} (recall Remark 
\ref{rem:change_center_3}). Then we can apply Corollary \ref{coro:quantitative_GN} and 
Proposition \ref{prop:location}, obtaining
\[
\begin{split}
\| \varphi_{\tilde{E}_{\varepsilon}} \|_{C^{1,\alpha}(\sphere^{N-1})} 
&\le C \left[\lambda^1(\tilde{E}_{\varepsilon}, \tilde{\Omega}_{\varepsilon}) - 
\lambda^1(B, \R^N) \right]^{(1-\alpha)/(4+N)} = C \left[\tilde\lambda_{\varepsilon} - \tilde\lambda_0 \right]^{(1-\alpha)/(4+N)}\\
&\le C \exp\left(-2\sqrt{\tilde\lambda_0\underline{m}}\cdot 2d^*\cdot\eps^{-1/N} 
\cdot \frac{1-\alpha}{4+N} \right),
\end{split}
\]
for every $0<\eps<\bar\eps$, and the conclusion follows.
\end{proof}

The proof of Theorem \ref{thm:quantitStabRNSharpPos} proceeds in several steps, that for the 
reader's convenience are briefly and conceptually listed as follows, each corresponding to 
a section below:
\begin{description}
    \item[{Step 0:}] we set the main notation and the basic properties of deformation paths
    (Section \ref{section:DeformationPathsProperties});  
    \item[{Step 1:}] we prove differentiability of eigenfunction and eigenvalue via the implicit function theorem (Section \ref{subsect:DifferentiabilitySharpPos});
    \item[{Step 2:}] for $C^{1, 1}$-nearly spherical domains, we carry out the expansion of the eigenvalue up to the second order (Section \ref{sec:expansion});    
    \item[{Step 3:}] for $C^{1, 1}$-nearly spherical domains, we prove coercivity of the second order shape derivative of the eigenvalue, computed at the origin, with respect to the 
    $L^2(\sphere^{N-1})$ norm of the deformation (Section \ref{sec:coercivity}); this is the point where the constraint on the barycenter plays a crucial role; 
    \item[{Step 4:}] for $C^{1, 1}$-nearly spherical domains, we prove continuity of the remainder with respect to the $L^2(\sphere^{N-1})$ norm of the deformation and a modulus of continuity depending on the $C^{1, 1}(\sphere^{N-1})$ error of the deformation (Section \ref{sec:cont_of_remaind}).
\end{description}

Once the above steps are performed, the proof of Theorem \ref{thm:quantitStabRNSharpPos} is 
straightforward and it is concluded at the end of Section \ref{sec:cont_of_remaind}.

As a final remark, we notice that Theorem \ref{thm:quantitStabRNSharpPos} is sharp, in the 
sense of equation \eqref{eq:sharp_indeed}.
\begin{lemma}\label{lem:sharp_indeed}
Under the assumptions of  there exists a constant $C'>0$ such that
\[
\lambda^1(\mathcal{A}, \R^N) - \lambda^1(B, \R^N) \le C' \| \varphi_{\mathcal{A}} \|^2_{L^2(\sphere^{N-1})}
\]
\end{lemma}
\begin{proof}
The proof is a direct computation, bounding $\lambda^1(\mathcal{A}, \R^N)$ from above with the associated Rayleigh quotient of $w$.

First, since $\Lcal(\Acal)=\Lcal(B)=1$, we have
\begin{equation}\label{eq:shi1}
0=\int_B d\bx - \int_\Acal d\bx = \int_{\sphere^{N-1} }\int_{r_0+\varphi_\Acal(\theta)}^{r_0} r^{N-1} dr\,d\sigma = \frac{1}{N}\int_{\sphere^{N-1} }\left[r_0^N-
\left(r_0+\varphi_\Acal(\theta)\right)^{N}\right]d\sigma
\end{equation}
which yields
\begin{equation}\label{eq:shi2}
\int_{\sphere^{N-1} }\varphi_\Acal(\theta) = - \frac{N-1}{2r_0} \| \varphi_{\mathcal{A}} \|^2_{L^2(\sphere^{N-1})} \left(1+R(\varphi_\Acal)\right),
\qquad\text{where }|R(\varphi_\Acal)|\le C_\delta
\end{equation}
(indeed, $R(\varphi_\Acal)=\sum_{n=3}^N\frac{2(N-2)!}{n!(N-n)!r_0^{n-2}}\int_{\sphere^{N-1} }\varphi_\Acal^{n-2}d\sigma$ and, for instance, $C_\delta=\sum_{n=3}^N\frac{2N(N-2)!\delta^{n-2}}{n!(N-n)!r_0^{n-2}}$).

Using the properties of $w$ and writing $m_\Acal = \overline{m}\ind{\Acal} - \underline{m}
\ind{\R^N\setminus\Acal} $ (and the same with $B$ instead of $\Acal$) we have 
\[
\begin{split}
\lambda^1(\mathcal{A}, \R^N) - \lambda^1(B, \R^N) &\le 
\frac{\int_{\R^N}|\nabla w|^2}{\int_{\R^N}m_\Acal w^2} - \frac{\int_{\R^N}|\nabla w|^2}{\int_{\R^N}m_B w^2}\\
&= \frac{\lambda^1(B, \R^N)}{\int_{\R^N}m_\Acal w^2}\int_{\R^N}(m_B - m_\Acal) w^2.
\end{split}
\]
Now, $\int_{\R^N}m_\Acal w^2 \ge \int_{\R^N}m_{B_{r_0-\delta}}w^2$, which is a positive constant independent of $\Acal$; on the other hand, since $w^2(r)\le A_\delta - B_\delta r$ 
for $r_0-\delta < r < r_0 + \delta$, with $A_\delta, B_\delta$ positive,
\[
\begin{split}
\int_{\R^N}(m_B - m_\Acal) w^2 &= (\overline{m}+\underline{m})\int_{\sphere^{N-1} }\int_{r_0+\varphi_\Acal(\theta)}^{r_0} r^{N-1}w^2(r) dr\,d\sigma\\ &\le -B_\delta (\overline{m}+\underline{m})\int_{\sphere^{N-1} }\int_{r_0+\varphi_\Acal(\theta)}^{r_0} r^{N} dr\,d\sigma
\end{split}
\]
where the term with $A_\delta$ cancels because of \eqref{eq:shi1}. Finally, we can evaluate the last integral and use \eqref{eq:shi2} to write
\[
\begin{split}
\int_{\sphere^{N-1} }\int_{r_0+\varphi_\Acal(\theta)}^{r_0} r^{N} dr\,d\sigma  
&= (N+1)r_0^{N+1} \left(-\int_{\sphere^{N-1} }\varphi_\Acal(\theta) - \frac{N}{2r_0} \| \varphi_{\mathcal{A}} \|^2_{L^2(\sphere^{N-1})} \left(1+R'(\varphi_\Acal)\right)\right) \\
&= - \frac{(N+1) r_0^N}{2} \| \varphi_{\mathcal{A}} \|^2_{L^2(\sphere^{N-1})} \left(1+R''(\varphi_\Acal)\right)\ge - C \| \varphi_{\mathcal{A}} \|^2_{L^2(\sphere^{N-1})},
\end{split}
\]
with $C>0$ independent on $\Acal$, and the lemma follows.
\end{proof}

\subsection{Deformation paths and their properties}\label{section:DeformationPathsProperties}

Before proceeding with the proof, we need to specify in what sense we intend the differentiation with respect to the shape, for instance of the eigenvalue. Analogous considerations hold for the associated eigenfunction, and for all the other quantities depending on a shape.

A standard procedure is to reduce to a one-dimensional problem, in a way we briefly recall. Consider a reference shape $\mathcal{A}_0 \in \R^N$, and a (small) perturbation $\mathcal{A} \in \R^N$. Suppose the existence of a one-parameter family of smooth displacement fields $\Phi(t, \mathbf{x}): [0, 1] \times \R^N \to \R^N$ which are diffeomorphisms for any fixed $t$, satisfying
\begin{equation}\label{eqn:CondFamDiffSharpIndef}
\Phi(0, \mathbf{x}) = \mathbf{I} \quad \text{and} \quad  \Phi(1, \mathcal{A}_0) = \mathcal{A} .
\end{equation}

At this point, one can introduce the parametric family
\[
\lambda(t) \coloneqq \lambda^1(\Phi(t, \mathcal{A}_0), \R^N) .
\]
Such quantity, once the reference shape has been fixed, depends only on the parameter $t$. Hence, one defines all the desired quantities exactly as in the case of functions of one variable. For instance, by first derivative of the eigenvalue, we mean the quantity
\[
\dot{\lambda}(t) \coloneqq \frac{d}{dt} \lambda(t) .
\]
Of course, one needs to check in advance that the desired quantities are well defined. In our case, this will be the subject of Section \ref{subsect:DifferentiabilitySharpPos}.

The existence of a family of diffeomorphisms satisfying \eqref{eqn:CondFamDiffSharpIndef} together with additional properties, for a general class of problems, has been the object of many works in the literature. An example is \cite{Dambrine2002:ShapeHessian}, where the vector fields $\partial_t \Phi(t, \mathbf{x})$ are required to be divergence free in a neighborhood of $\partial \mathcal{A}_0$.

As remarked in \cite[Appendix A]{Brasco:SharpFaberKrahn}, when the reference configuration $\mathcal{A}_0$ coincides with 
a ball and when nearly spherical perturbations are considered, as for our problem, the construction of a family 
$\Phi(t, \mathbf{x})$ can be made explicit. In radial coordinates, it is sufficient to consider 
\begin{equation}\label{eqn:PolarDisplacFamilySharpIndef}
\Phi_\Acal(t, \rho \, \theta) \coloneqq \left[ \rho^N + t \left( \, ( r_0+\varphi_{\mathcal{A}}(\theta) )^{N} - r_0^N \right) \right]^{1/N} \, \theta \quad \text{ for } \theta \in \sphere^{N-1}, \, 3\frac{r_0}{4} < \rho < 5 \frac{r_0}{4}.
\end{equation}
Then, one can extend \eqref{eqn:PolarDisplacFamilySharpIndef} to the whole $\R^N$, so that it also satisfies some additional properties. More precisely, similarly as in  \cite[Appendix A]{Brasco:SharpFaberKrahn}, we have the following.
\begin{lemma}\label{lemma:DivFreeSpherDeformSharpIndef}
There exists $\delta>0$ sufficiently small and a modulus of continuity $\eta$ such that, for any set $\mathcal{A} \subset \R^N$ nearly spherical of class $C^{1, 1}$ with
\begin{equation}\label{eqn:LemmaSherDefConditSharpIndef}
\| \varphi_{\mathcal{A}} \|_{C^{1, 1}(\sphere^{N-1})} < \delta  \quad  \text{and} \quad  \mathcal{L}(\mathcal{A}) = \mathcal{L}(B)=1 ,
\end{equation}
there exists a family $\Phi_{\mathcal{A}}(t, \rho \, \theta)$ extending \eqref{eqn:PolarDisplacFamilySharpIndef} to the whole $\R^N$, such that
\begin{itemize}
    \item [(i)] The map $t \in [0, 1] \mapsto \Phi_{\mathcal{A}}(t, \mathbf{x}) - \mathbf{I} \in W^{2, \infty}(\R^N)$ is of class $C^{\infty}$.
    \item [(ii)] For $3r_0/4<|\bx|<5r_0/4$, 
    \[
    \partial_t \Phi_{\mathcal{A}}(t, \mathbf{x}) = \mathbf{X}_{\mathcal{A}}( \Phi_{\mathcal{A}}(t, \mathbf{x}) ), \qquad \diverg(\mathbf{X}_{\mathcal{A}}) = 0,
    \]
    where
    \begin{equation}\label{enq:DefXSharpIndef}
    \mathbf{X}_{\mathcal{A}}(\rho, \theta) \coloneqq \frac{( r_0+\varphi_{\mathcal{A}}(\theta) )^N - r_0^N}{N \, \rho^{N-1}} \, \theta .
    \end{equation}
    \item [(iii)] $\mathcal{L}(\Phi_{\mathcal{A}}(t, B))=1$ for all $t \in [0, 1]$.
    \item [(iv)] $\| \Phi_{\mathcal{A}} - \mathbf{I} \|_{C^{1, 1}(\R^N)} \le \eta\left( \| \varphi_{\mathcal{A}} \|_{C^{1, 1}(\sphere^{N-1})} \right)$ uniformly for $t \in [0, 1]$.
    \item[(v)] $\mathbf{X}_{\mathcal{A}} \circ \Phi_{\mathcal{A}} - \mathbf{X}_{\mathcal{A}} = O \left( \| \varphi_{\mathcal{A}} \|_{C^{1, 1}(\sphere^{N-1})} \right) \mathbf{X}_{\mathcal{A}}$ on $\partial B$, uniformly for $t \in [0, 1]$.
    \item[(vi)] $\left\| \mathbf{X}_{\mathcal{A}}(r_0,\cdot) \cdot \frac{\mathbf{x}}{|\bx|} - \varphi_{\mathcal{A}}  \right\|_{L^2(\sphere^{N-1})} \le \eta \left( \| \varphi_{\mathcal{A}} \|_{L^{\infty}(\sphere^{N-1})} \right) \| \varphi_{\mathcal{A}}  \|_{L^2(\sphere^{N-1})}$.
\end{itemize}
\end{lemma}
\begin{remark}
In Lemma \ref{lemma:DivFreeSpherDeformSharpIndef} (i) the derivatives as $t=0$ and $t=1$ are understood from the right and from the left, respectively. Moreover, in Lemma \ref{lemma:DivFreeSpherDeformSharpIndef} (iv) and (v) analogous estimates continue to hold if the $C^{1, 1}(\partial B)$ norm is substituted by the $C^{1, \alpha}(\partial B)$ norm, for any $0 < \alpha < 1$.
\end{remark}
\begin{proof}
Consider a smooth cutoff function $h$, such that $h \equiv 1$ for $3r_0/4 \le \rho \le 5 r_0/4$ and $h \equiv 0$ for $r_0/2 \le \rho \le 3 r_0/2$. We define
\begin{equation}\label{eqn:DeformExprSharpIndef}
\Phi_{\mathcal{A}}(t, \rho \, \theta) \coloneqq \left[ \rho^N + t \, h(\rho) \left( \, ( r_0+\varphi_{\mathcal{A}}(\theta) )^{N} - r_0^N \right) \right]^{1/N} \, \theta
\end{equation}
in the whole $\R^N$. Clearly, \eqref{eqn:DeformExprSharpIndef} is an extension of \eqref{eqn:PolarDisplacFamilySharpIndef}.

Let us prove point $(i)$. Thanks to the restriction on $\| \varphi_{\mathcal{A}} \|_{C^{1, 1}(\sphere^{N-1})}$ required in \eqref{eqn:LemmaSherDefConditSharpIndef}, for $\delta > 0$ sufficiently small it holds that the quantities $\Phi_{\mathcal{A}}(t, \mathbf{x}) - \mathbf{I}$, $ D_{\mathbf{x}} \Phi_{\mathcal{A}}(t, \mathbf{x}) - \mathbf{I}$ and $ D_{\mathbf{x}}^2 \Phi_{\mathcal{A}}(t, \mathbf{x})$ are $C^{\infty}$ in time in the classical sense, for almost every $\mathbf{x}=(\rho, \theta) \in \R^N$ fixed. By a simple computation, all time derivatives are bounded in $L^{\infty}(\R^N)$, uniformly for $t \in [0, 1]$. This is sufficient to prove point $(i)$.

Point $(ii)$ follows immediately by direct computation.

Point $(iii)$ follows from the fact that, by \cite[Corollaire 5.2.8]{HenrotPierre:ShapeOptimiz}, it holds that
\begin{equation}\label{eqn:ConservedVolumeSharpIndef}
\frac{d}{dt} \mathcal{L}(\Phi_{\mathcal{A}}(t, B)) = \frac{d}{dt} \int_{\Phi_{\mathcal{A}}(t, B)} 1 = \int_{\partial  \Phi_{\mathcal{A}}(t, B)} \mathbf{X}_{\mathcal{A}} \cdot \mathbf{n} ,
\end{equation}
where in the last passage we have used points $(i)$ and $(ii)$. Now, using \eqref{enq:DefXSharpIndef}, the last integral in \eqref{eqn:ConservedVolumeSharpIndef} vanishes thanks to the condition $\mathcal{L}(\mathcal{A}) = \mathcal{L}(B)$ in \eqref{eqn:LemmaSherDefConditSharpIndef}. Hence, point $(iii)$ follows.

Point $(iv)$ follows easily by direct computation.

Point $(v)$ also follows by direct computation, using point $(iv)$.

Point $(vi)$ follows expanding the expression \eqref{enq:DefXSharpIndef} for $\mathbf{X}_{\mathcal{A}}$ in terms of $\varphi_{\mathcal{A}}$. Notice that, since $\| \varphi_{\mathcal{A}}  \|_{L^2(\sphere^{N-1})} \le 2 \| \mathbf{X}_{\mathcal{A}}(r_0,\cdot) \cdot \frac{\mathbf{x}}{|\bx|} \|_{L^2(\sphere^{N-1})} $ for $\delta$ sufficiently small, up to redefining $\eta$ it also holds that
\begin{equation}\label{eqn:EstXvarphiDiffL2SharpIndef}
    \left\| \mathbf{X}_{\mathcal{A}}(r_0,\cdot) \cdot \frac{\mathbf{x}}{|\bx|} - \varphi_{\mathcal{A}}  \right\|_{L^2(\sphere^{N-1})} \le \eta \left( \| \varphi_{\mathcal{A}} \|_{L^{\infty}(\sphere^{N-1})} \right) \left\| \mathbf{X}_{\mathcal{A}}(r_0,\cdot) \cdot \frac{\mathbf{x}}{|\bx|} \right\|_{L^2(\sphere^{N-1})} .
\end{equation}
\end{proof}
\begin{remark}[on notation]\label{remark:PathNotationSharpIndef}
In the following, we will denote with $\Phi_{\mathcal{A}}(t, \rho \theta)$ or equivalently $\Phi_{\mathcal{A}}(t, \mathbf{x})$ the family of diffeomorphic deformations of the form \eqref{eqn:PolarDisplacFamilySharpIndef} already extended to the whole $\R^N$, in order to satisfy the conditions of Lemma \ref{lemma:DivFreeSpherDeformSharpIndef}.
\end{remark}

\subsection{Differentiability}\label{subsect:DifferentiabilitySharpPos}

Consider $\Theta \in W^{1, \infty}(\R^N, \R^N)$, with $\| \Theta \|_{W^{1, \infty}(\R^N, \R^N)}$ small enough so that $\mathbf{I} + \Theta$ is a diffeomorphism of $\R^N$ in itself. We introduce the domain $D_{\Theta}$ defined as the image of $B$ via $\mathbf{I} + \Theta$, the correspondent weight 
\[
\tilde{m}_{\Theta} \coloneqq \overline{m}_{\mathcal{X}_{D_{\Theta}}} -\underline{m}_{\mathcal{X}_{ \R^N \setminus D_{\Theta} } }
\]
and the eigenfunction $u_{\Theta}$, solution in $H^1(\R^N)$ of
\begin{equation}\label{eqn:DiffPbThetaRNSharpPos}
-\Delta u_{\Theta} = \lambda^1(D_{\Theta}, \R^N) \tilde{m}_{\Theta} u_{\Theta} \quad \text{in } \R^N ,
\end{equation}
with the normalization condition
\[
\int_{\R^N} \tilde{m}_{\Theta} \, u_{\Theta}^2 = 1 .
\]
A common practice to study the differentiability of $u_{\Theta}$ with respect to $\Theta$, is to bring back the problems for different vector fields $\Theta$ on the same domain $B$ and study first the properties of such problem, which are usually easier to tackle. Then, one tries to transfer the obtained information to the function $u_{\Theta}$. This can be done introducing the functions
\[
v_{\Theta} \coloneqq u_{\Theta} \circ (\mathbf{I} + \Theta).
\]
We remark that, thanks to the regularity required to $\Theta$, it holds that $v_{\Theta} \in H^1(\R^N)$. Moreover, performing a change of variables in the weak formulation of \eqref{eqn:DiffPbThetaRNSharpPos}, it is easy to see that the functions $v_{\Theta}$ solve, in $H^1(\R^N)$
\begin{equation}\label{eqn:DiffPbB1ThetaRNSharpPos}
-\diverg( A(\Theta) v_{\Theta} ) = \lambda^1(D_{\Theta}, \R^N) J(\Theta) \tilde{m}_0 v_{\Theta} \quad \text{in } \R^N ,
\end{equation}
where $J(\Theta)$ denotes the jacobian of the diffeomorphism $\mathbf{I}+\Theta$, i.e.
\[
J(\Theta) \coloneqq \det (\mathbf{I}+D\Theta)
\]
while $A(\Theta)$ is defined as
\[
A(\Theta) \coloneqq J(\Theta) (I+D\Theta)^{-1} (I+D\Theta)^{-T} .
\]

The main result of this section concerns a local differentiability of the function 
$v_{\Theta}$, and of the corresponding eigenvalue, with respect to $\Theta$. To 
prove such result we need a preliminary lemma.
\begin{lemma}\label{lem:LM}
For every $h\in L^2(\R^N)$ such that $\int_{\R^N} w h=0$ there exists a unique $y\in H^2(\R^N)$, $\int_{\R^N}\tilde m_0 w y=0$, such that 
\[
- \Delta z - \tilde{\lambda}_0 \tilde{m}_0 z = h
\quad\iff\quad
z=y+cw,\ c\in\R.
\]
On the other hand, if $\int_{\R^N} w h\neq 0$ the above equation has no solution.
\end{lemma}
\begin{proof}
The necessary condition on $h$ for the existence of solutions follows by testing the equation with $w$ and using the corresponding equation.

To show existence, let us introduce the bilinear form on $H^1(\R^N)$:
\[
a(y,v)\coloneqq\int_{\R^N} \nabla y \cdot \nabla v - \tilde{\lambda}_0 \tilde{m}_0 y v.
\]
By the variational characterization of $\tilde \lambda_0>0$ and $w$ we have that   
\begin{equation}\label{eq:varyy1}
a(v,v) \le 0 
\quad\iff\quad
v=cw,\ c\in\R
\end{equation}
(and hence $a(v,v)=0$). Moreover, considering the Hilbert space
\[
V\coloneqq \left\{v\in H^1(\R^N):\int_{\R^N}\tilde m_0 w v=0\right\},
\qquad
\|v\|_V^2 \coloneqq \|\nabla v\|_{L^2(\R^N)}^2 + \tilde \lambda_0 \underline{m}
\| v\|_{L^2(\R^N)}^2,
\]
we have that
\begin{equation}\label{eq:varyy2}
a(v,v) = \|v\|^2_V - \tilde\lambda_0(\overline{m} + \underline{m}) \|v\|_{L^2(B)}^2;
\end{equation}

Then we can apply the Lax-Milgram theorem to the variational problem: 
\begin{equation}\label{eq:var_V}
\text{find $y\in V$ such that }\qquad 
a(y,v) =\int_{\R^N} hv\qquad\text{for every }v\in V.
\end{equation}
Indeed, (bi)linearity and continuity are straightforward, while the coercivity 
of $a$ follows by contradiction: assume the existence of a sequence 
$(v_n)_n\subset V$ such that
\[
\|v_n\|_V=1,\qquad a(v_n,v_n)\le \frac{1}{n}.
\]
Then, up to subsequences, $v_n$ converges to some $\bar v$, weakly in $V$ and strongly in $L^2_{\loc}(\R^N)$, as $n\to+\infty$. In particular, exploiting the lower weak semicontinuity of the norm in \eqref{eq:varyy2} we infer $a(\bar v,\bar v)\le0$, and, since $\bar v\in V$, \eqref{eq:varyy1} yields $\bar v \equiv 0$; on the other hand, 
again by \eqref{eq:varyy2},
\[
\tilde\lambda_0(\overline{m} + \underline{m}) \|\bar v\|_{L^2(B)}^2=
\tilde\lambda_0(\overline{m} + \underline{m}) \|v_n\|_{L^2(B)}^2 + o(1)
= \|v_n\|_V^2-a(v_n,v_n) + o(1) \ge 1 -\frac{1}{n}+ o(1)
\]
for $n$ large, contradiction.

Then Lax-Milgram theorem applies, yielding the existence and uniqueness of 
$y\in V$ satisfying \eqref{eq:var_V}. Taking $v\in H^1(\R^N)$ we have that 
$v-t_v w\in V$, where $t_v=\int_{\R^N} \tilde m_0 wv / \int_{\R^N}
 \tilde m_0 w^2$. Substituting in \eqref{eq:var_V} we have
\[
a(y,v) = \int_{\R^N} hv - t_v \left[a(y,w) - \int_{\R^N} hw \right]
= \int_{\R^N} hv
\qquad \text{for every }v\in H^1(\R^N),
\]
where we used the equation of $w$ and the fact that $\int_{\R^N} w h=0$ 
by assumption. Then by elliptic regularity $y\in H^2(\R^N)$, and $z=y+cw$, 
$c\in\R$, satisfies the equation
\[
- \Delta z - \tilde{\lambda}_0 \tilde{m}_0 z = h.
\]

Finally, consider $z_1,z_2\in H^1(\R^N)$ solutions of the above 
equation. Then subtracting we have $a(z_1-z_2,z_1-z_2)=0$ and we can conclude the proof using again \eqref{eq:varyy1}.
\end{proof}
\begin{proposition}\label{prop:W2inftyRegulSharpPos}
The function
\[
\Theta \in W^{2, \infty}(\R^N, \R^N) \mapsto (\lambda^1(D_{\Theta}, \R^N),  v_{\Theta}) \in \R \times H^2(\R^N)
\]
is $C^{\infty}$ in a neighborhood of $\Theta = \mathbf{0}$.
\end{proposition}
\begin{proof}
The proof follows from the implicit function theorem, similarly as in  \cite{Dambrine2011:shapeSensitivity}.

Let us define the map $F: W^{2, \infty}(\R^N, \R^N) \times \R \times H^2(\R^N) \to L^2(\R^N) \times \R$ as follows:
\[
F(\Theta, \lambda, v) \coloneqq 
\left( 
-\diverg( A(\Theta) \nabla v ) - \lambda \tilde{m}_0 v J(\Theta), 
\int_{\R^N} \tilde{m}_0 \, v^2 \, J(\Theta) - 1
 \right) .
\]
Then there exists a neighborhood $\Ucal \subset W^{2, \infty}(\R^N, \R^N)$ of the origin such that 
\begin{equation}\label{eqn:CInftyFSharpPos}
F(\Theta, \lambda, v) \in C^{\infty}(\Ucal \times \R \times H^2(\R^N)) .
\end{equation}
Indeed, it is a standard result (reasoning for instance as in \cite[Théorème 5.3.2]{HenrotPierre:ShapeOptimiz} and its proof) that the maps
\[
(\Theta, v) \in \Ucal \times H^2(\R^N) \mapsto \diverg( A(\Theta) \nabla v ) \in L^2(\R^N) ,
\]
\[
\Theta \in W^{2, \infty}(\R^N, \R^N) \mapsto J(\Theta) \in \R
\]
are of class $C^{\infty}$, in an appropriate $\Ucal$. Thanks to the regularity of $J(\Theta)$, it is immediate to deduce the $C^{\infty}$ regularity property also for
\[
(\Theta, \lambda, v) \in \Ucal  \times \R \times H^2(\R^N) \mapsto \lambda \tilde{m}_0 v J(\Theta) \in L^2(\R^N) ,
\]
\[
(\Theta, v) \in \Ucal  \times H^2(\R^N) \mapsto \int_{\R^N} \tilde{m}_0 \, v^2 \, J(\Theta) \in \R .
\]
Hence, \eqref{eqn:CInftyFSharpPos} is proved.

Next, recalling that $J(\mathbf{0})=1$ and $A(\mathbf{0}) = \mathbf{I}$, we easily 
infer that, for any $\mu \in \R$ and $z \in H^2(\R^N)$,
\[
D_{\lambda, v}(\mathbf{0}, \tilde{\lambda}_0, w)[\mu, z] = \left( 
- \Delta z - \mu \tilde{m}_0 w - \tilde{\lambda}_0 \tilde{m}_0 z , 
\int_{\R^N} 2\tilde{m}_0 \, z \, w 
 \right).
\]

In order to apply the implicit function theorem, we are left to prove that $D_{\lambda, v}(\mathbf{0}, \tilde{\lambda}_0, w)$ is a diffeomorphism. Actually, proving its invertibility is sufficient since it is continuous and linear, hence the continuity of its inverse would follow from the inverse mapping theorem in Banach spaces.

To prove invertibility, consider a couple $(k, p) \in L^2(\R^N) \times \R$. We have to show existence and uniqueness of the pair 
$(\mu, z) \in \R \times H^2(\R^N)$ solving
\[
\begin{cases}
    - \Delta z - \tilde{\lambda}_0 \tilde{m}_0 z = k + \mu \tilde{m}_0 w , \\ 
    \int_{\R^N} 2\tilde{m}_0 \, z \, w = p . 
\end{cases}
\]
The first equation is solvable if and only if the condition on the right 
hand side appearing in Lemma \ref{lem:LM} is satisfied, i.e.\ 
\[
\mu = - \frac{\int_{\R^N} k w}{\int_{\R^N} \tilde m_0 w^2}.
\]
with infinitely many solutions given by
\[
z = y + c w,\qquad c\in\R,
\]
where $y$ is uniquely determined in such a way that 
$\int_{\R^N} \tilde m_0 w y=0$. Substituting in the second equation we have that 
$c$ is uniquely determined by
\[
c = \frac{p}{2\int_{\R^N} \tilde m_0 w^2},
\]
concluding the proof.
\end{proof}

Our aim now is to exploit the differentiability of $v_{\Theta}$ to prove a corresponding result for $u_{\Theta}$. Actually, since as recalled in Section \ref{section:DeformationPathsProperties} we will exploit differentiation along paths, we derive the result directly for $u_{\Theta}$ restricted on such paths.

Let $\mathcal{A}$ be a nearly spherical set of class $C^{1, 1}$, with $\| \varphi_{\mathcal{A}} \|_{C^{1, 1}(\sphere^{N-1})}$ sufficiently small. Consider the path $\Phi_{\mathcal{A}}(t, \mathbf{x})$ associated to $\mathcal{A}$ (see also Remark \ref{remark:PathNotationSharpIndef}) and introduce the parametric family $\mathcal{A}_{t} \coloneqq \Phi_{\mathcal{A}}(t, B)$, the weights $m_{t} \coloneqq \overline{m}\ind{\Acal_t}-\underline{m}\ind{\R^N\setminus \Acal_t}$ and the functions $u_t \in H^1_0(\R^N)$ for any fixed $t \in [0, 1]$, solutions of
\begin{equation}\label{eqn:FlowDiffProblemSharpIndef}
-\Delta u_t = \lambda^1(\mathcal{A}_t, \R^N) m_t u_t \qquad \text{in } \R^N
\end{equation}
As usual, the principal eigenvalue $\lambda^1(\mathcal{A}_t, \R^N)$ in \eqref{eqn:FlowDiffProblemSharpIndef} has a unique eigenfunction $u_t$, up to normalization, see e.g.\ 
\cite[Lemma 2.6]{Verzini:Neumann}.

Then, from Proposition \ref{prop:W2inftyRegulSharpPos} we can deduce the following.
\begin{lemma}\label{lem:uC1TImeSharpIndef}
There exists $\delta>0$ sufficiently small such that, for any nearly spherical set $\mathcal{A}$ of class $C^{1, 1}$ with $\| \varphi_{\mathcal{A}} \|_{C^{1, 1}(\sphere^{N-1})} < \delta$, the map $t \in [0, 1] \mapsto u_t \in H^1(\R^N)$ is of class $C^1$.
\end{lemma}
\begin{proof}
Let us denote $v_{t}$ the map $v_t \coloneqq v_{\Phi_{\mathcal{A}}(t) - \mathbf{I}}$. Thanks to Lemma \ref{lemma:DivFreeSpherDeformSharpIndef} and Proposition \ref{prop:W2inftyRegulSharpPos}
, by composition it holds that $t\in[0,1] \mapsto v_t\in H^2(\R^N)$ is of class $C^{\infty}$.

Now, the function $u_t$ can be written as $u_t = v_t \circ \Phi_{\mathcal{A}}(t, \mathbf{x})^{-1}$. Since the map $\Theta \in W^{2, \infty} \mapsto (\mathbf{I}+\Theta)^{-1} \in W^{1, \infty}$ is of class $C^1$ in a neighborhood of $\Theta = \mathbf{0}$, by composition the map $t\mapsto \Phi_{\mathcal{A}}(t, \mathbf{x})^{-1}\in  W^{1, \infty}$ is $C^1$ in $[0, 1]$ (see e.g.\ 
\cite[Sec. 5.2.2]{HenrotPierre:ShapeOptimiz}). Summing up, we deduce that
\begin{equation}\label{eq:HP1}
t \mapsto (v_t, \Phi_{\mathcal{A}}(t, \mathbf{x})^{-1} - \mathbf{I}) \in H^2(\R^N) \times W^{1, \infty}
\qquad\text{ is of class $C^1$.}
\end{equation}
On the other hand, by \cite[Lemme 5.3.9]{HenrotPierre:ShapeOptimiz}, in a sufficiently small neighborhood $\Ucal  \subset W^{1, \infty}$ of the origin, the map 
\begin{equation}\label{eq:HP2}
(g, \Theta) \in H^2(\R^N) \times \Ucal  \mapsto g \circ (\mathbf{I}+\Theta) \in H^1(\R^N)\qquad\text{ is of class $C^1$.}
\end{equation}
Hence, using Lemma \ref{lemma:DivFreeSpherDeformSharpIndef} (iv), for $\delta$ sufficiently small we can conclude the proof by composition of the maps in \eqref{eq:HP1} and \eqref{eq:HP2}.
\end{proof}

\subsection{The expansion}\label{sec:expansion}

In this section we perform the expansion of the eigenvalue up to the second order. We will use extensively the Hadamard formula, see e.g. \cite[Corollaire 5.2.8]{HenrotPierre:ShapeOptimiz}. 
In particular, let
\[
f(x, t) \in C^1\left([0, T); L^1(\R^N)\right) \cap C^0\left([0, T); W^{1, 1}(\R^N)  \right) ,
\]
then, recalling the definition and properties of $\mathbf{X}_{\mathcal{A}}$ 
(see Lemma \ref{lemma:DivFreeSpherDeformSharpIndef}),
\begin{align*}
   &\frac{d}{dt} \int_{\mathcal{A}_t} f(x, t) = \int_{\mathcal{A}_t} \partial_t f(x, t) + \int_{\mathcal{A}_t} \diverg \left( f(x, t) \mathbf{X}_{\mathcal{A}}(x) \right) ,\\
   &\frac{d}{dt} \int_{\R^N \setminus \mathcal{A}_t} f(x, t) = \int_{\R^N \setminus \mathcal{A}_t} \partial_t f(x, t) + \int_{\R^N \setminus \mathcal{A}_t} \diverg \left( f(x, t) \mathbf{X}_{\mathcal{A}}(x) \right).
\end{align*}

In the following we denote with 
$\lambda_t \coloneqq \lambda^1(\mathcal{A}_t, \R^N)$, which is of class 
$C^{\infty}$ in time by Proposition \ref{prop:W2inftyRegulSharpPos}, and with $u_t$ 
the corresponding eigenfunction as in \eqref{eqn:FlowDiffProblemSharpIndef} which 
is $C^1$ in time, with values in $H^1$, by Lemma \ref{lem:uC1TImeSharpIndef}. Moreover, we denote with $\dot\lambda_t$, $\dot{u}_t$ the time derivatives of 
$\lambda_t$ and of $u_t$, respectively.

We start with the first order expansion.
\begin{lemma}\label{eqn:1stOrderEigenExpSharpIndef}
There exists $\delta>0$ sufficiently small such that, for any nearly spherical set $\mathcal{A}$ of class $C^{1, 1}$ with $\| \varphi_{\mathcal{A}} \|_{C^{1, 1}(\sphere^{N-1})} < \delta$
\begin{equation}\label{enq:1stOrderEigenExpSharpIndef}
\dot{\lambda}_t =  - \lambda_t (\overline{m}+\underline{m}) \int_{\partial \mathcal{A}_t} u_t^2 \, \mathbf{X}_{\mathcal{A}} \cdot \mathbf{n}_t ,
\end{equation}
where $\mathbf{n}_t$ denotes the unitary outer normal to $\partial \mathcal{A}_t$. 
\end{lemma}
\begin{proof}
Choosing the normalization $\int_{\R^N} m_t u_t^2 = 1$ we have that
\begin{equation}\label{eqn:eqn1stOrderEigenExpSharpIndef}
    \dot{\lambda}_t  = \frac{d}{dt} \int_{\R^N} \vert \nabla u_t \vert^2 = 2 \int_{\R^N} \nabla \dot u_t \nabla u_t .
\end{equation}
Our aim is to rewrite the previous formula only in terms of $u_t$. Exploiting the chosen normalization and using the Hadamard formula, we have that
\begin{align}\label{eqn:eqn2stOrderEigenExpSharpIndef}
\begin{split}
    & 0 = \frac{d}{dt} \int_{\R^N} m_t u_t^2 = \frac{d}{dt} \left[ \overline{m} \int_{\mathcal{A}_t} u_t^2 - \underline{m} \int_{\mathcal{A}^c_t} u_t^2 \right] = \\
    & = 2 \int_{\R^N} m_t \dot u_t u_t + (\overline{m} + \underline{m}) \int_{\partial \mathcal{A}_t} u_t^2 \, \mathbf{X}_{\mathcal{A}} \cdot \mathbf{n}_t .
    \end{split}
\end{align}
Testing \eqref{eqn:FlowDiffProblemSharpIndef} with $\dot u_t$, we get
\begin{equation}\label{eqn:eqn3stOrderEigenExpSharpIndef}
    \int_{\R^N} \nabla \dot u_t \nabla u_t = \lambda_t \int_{\R^N} m_t \dot u_t u_t .
\end{equation}
Substituting \eqref{eqn:eqn3stOrderEigenExpSharpIndef} in \eqref{eqn:eqn1stOrderEigenExpSharpIndef} and the using \eqref{eqn:eqn2stOrderEigenExpSharpIndef}, we can conclude.
\end{proof}

Now we differentiate the eigenvalue once again, to obtain the following.
\begin{lemma}\label{lemma:2ndOrderEigenExpSharpIndef}
There exists $\delta>0$ sufficiently small such that, for any nearly spherical set $\mathcal{A}$ of class $C^{1, 1}$ with $\| \varphi_{\mathcal{A}} \|_{C^{1, 1}(\sphere^{N-1})} < \delta$
\begin{align}\label{eqn:2ndOrdertEigenExpSharpIndef}
\begin{split}
& \ddot{\lambda}_t  =  \lambda_t \left[ (\overline{m}+\underline{m}) \int_{\mathcal{A}_t} \diverg(u_t^2 \mathbf{X}_{\mathcal{A}}) \right]^2 + \\
& - 2 \lambda_t (\overline{m}+\underline{m}) \left[ \int_{\partial \mathcal{A}_t} u_t \left(\nabla u_t \cdot  \mathbf{X}_{\mathcal{A}}\right) \left(\mathbf{X}_{\mathcal{A}} \cdot \mathbf{n}_t\right)  + \int_{\partial \mathcal{A}_t} \dot u_t u_t \mathbf{X}_{\mathcal{A}} \cdot \mathbf{n}_t \right]  .
\end{split}
\end{align}
where $\mathbf{n}_t$ denotes the unitary outer normal to $\partial \mathcal{A}_t$. In particular,
\begin{equation}\label{eqn:2ndOrder0EigenExpSharpIndef}
\ddot{\lambda}_t \vert_{t = 0} = 
2 \tilde \lambda_0 (\overline{m}+\underline{m}) w\vert_{\partial B} \left[ \left\vert \frac{\partial w}{\partial \rho} \right\vert_{\partial B} \int_{\partial B} \vert \mathbf{X}_{\mathcal{A}} \vert^2 - \int_{\partial B} \dot w \mathbf{X}_{\mathcal{A}} \cdot \mathbf{n}_0 \right]  .
\end{equation}
\end{lemma}
\begin{proof}
We rewrite \eqref{enq:1stOrderEigenExpSharpIndef} as 
\[
\dot{\lambda}_t =  - \lambda_t (\overline{m}+\underline{m}) \int_{\mathcal{A}_t} \diverg \left( u_t^2 \, \mathbf{X}_{\mathcal{A}} \right) ,
\]
and differentiating with respect to $t$ using the Hadamard formula, we get
\begin{align*}
& \ddot \lambda_t = - \dot \lambda_t (\overline{m}+\underline{m}) \int_{\mathcal{A}_t} \diverg \left( u_t^2 \, \mathbf{X}_{\mathcal{A}} \right) + \\
& - 2 \lambda_t (\overline{m}+\underline{m}) \int_{\mathcal{A}_t} \diverg \left( \dot u_t u_t \, \mathbf{X}_{\mathcal{A}} \right) - \lambda_t (\overline{m}+\underline{m}) \int_{\partial \mathcal{A}_t} \diverg \left( u_t^2 \, \mathbf{X}_{\mathcal{A}} \right) \mathbf{X}_{\mathcal{A}} \cdot \mathbf{n}_t = \\
& = - \dot \lambda_t (\overline{m}+\underline{m}) \int_{\mathcal{A}_t} \diverg \left( u_t^2 \, \mathbf{X}_{\mathcal{A}} \right) - 2 \lambda_t (\overline{m}+\underline{m}) \int_{\mathcal{A}_t} \diverg \left( \dot u_t u_t \, \mathbf{X}_{\mathcal{A}} \right) + \\
& - 2 \lambda_t (\overline{m}+\underline{m}) \int_{\partial \mathcal{A}_t} u_t \left(\nabla u_t \cdot  \mathbf{X}_{\mathcal{A}}\right)\left( \mathbf{X}_{\mathcal{A}} \cdot \mathbf{n}_t\right) - \lambda_t (\overline{m}+\underline{m}) \int_{\partial \mathcal{A}_t} u_t^2 \, \diverg \left( \mathbf{X}_{\mathcal{A}} \right) \mathbf{X}_{\mathcal{A}} \cdot \mathbf{n}_t .
\end{align*}
Substituting $\dot \lambda_t$ with \eqref{enq:1stOrderEigenExpSharpIndef} and noticing that the last term vanishes thanks to Lemma \ref{lemma:DivFreeSpherDeformSharpIndef} (ii), gives \eqref{eqn:2ndOrdertEigenExpSharpIndef}.

To obtain \eqref{eqn:2ndOrder0EigenExpSharpIndef} it is sufficient to notice that $u_0 = w$, $\mathbf{n}_0 = \bx/|\bx|$ and that the first term vanishes since $\partial B$ is a level set for $w$ and thanks to Lemma \ref{lemma:DivFreeSpherDeformSharpIndef} (ii).
\end{proof}

To conclude this section we derive an equation for $\dot u_t$.
\begin{lemma}\label{lemma:eqnDotutSharpInder}
For any $t \in [0, 1]$, the function $\dot u_t$ solves, in $H^1(\R^N)$, the following problem:
\begin{equation}\label{eqn:eqnDotutTSharpInder}
\begin{cases}
-\Delta \dot u_t = \dot \lambda_t m_t u_t + \lambda_t m_t \dot u_t  & \text{in } \R^N\setminus \partial \mathcal{A}_t , \\
[\dot u_t] = 0, \quad [\nabla \dot u_t \cdot \mathbf{n}_t] =  \lambda_t (\overline{m}+\underline{m}) u_t \, \mathbf{X}_{\mathcal{A}} \cdot \mathbf{n}_t  & \text{on } \partial \mathcal{A}_t \\
2 \int_{\R^N} m_t \dot u_t u_t = - (\overline{m} + \underline{m}) \int_{\partial \mathcal{A}_t} u_t^2 \, \mathbf{X}_{\mathcal{A}} \cdot \mathbf{n}_t ,
\end{cases}
\end{equation}
and the jump is defined as in \eqref{eq:jumpop}.

In particular, for $t = 0$ it holds that $u_0=w$, $\mathbf{n}_0=\bx/|\bx|$ and
\begin{equation}\label{eqn:eqnDotut0SharpInder}
\begin{cases}
-\Delta \dot w = \tilde{\lambda}_0 \tilde{m}_0 \dot w  & \text{in } \R^N \setminus \partial B, \\
[\dot w] = 0, \quad [\partial_{\rho} \dot w] =  \tilde{\lambda}_0 (\overline{m}+\underline{m}) w \, \mathbf{X}_{\mathcal{A}} \cdot \mathbf{n}_0  & \text{on } \partial B \\
2 \int_{\R^N} \tilde m_0 \dot w w = - (\overline{m} + \underline{m}) \int_{\partial B} w^2 \, \mathbf{X}_{\mathcal{A}} \cdot \mathbf{n}_0 = 0 ,
\end{cases}
\end{equation}
with $[\partial_{\rho} \dot w](r_0\theta)=\partial_{\rho} \dot w(r_0^-\theta)-\partial_{\rho} 
\dot w(r_0^+\theta)$.
\end{lemma}
\begin{proof}
The proof is straightforward, writing the weak formulation of \eqref{eqn:FlowDiffProblemSharpIndef} for any test function $\psi \in C^{\infty}_0(\R^N)$, and differentiating the formulation using the Hadamard formula.
\end{proof}

\subsection{Coercivity}\label{sec:coercivity}

In this section we prove the following result.
\begin{proposition} \label{prop:CoercRemainderSharpIndef}
There exist $\delta>0$ sufficiently small and a positive constant $C>0$ depending only on 
$\overline{m}$, $\underline{m}$ and the dimension $N$ such that, for any nearly spherical set $\mathcal{A}$ of class $C^{1, 1}$ with $\| \varphi_{\mathcal{A}} \|_{C^{1, 1}(\sphere^{N-1})} < \delta$, $\Lcal(\Acal)=1$  and $\bari(\mathcal{A}) = \mathbf{0}$, it holds 
\begin{equation}\label{eqn:CoercivityEstimSHarpIndef}
\left.\ddot \lambda_t\right|_{t=0}  \ge C\left \|  \mathbf{X}_{\mathcal{A}}(r_0,\cdot) \cdot \mathbf{n}_0\right \|_{L^2(\sphere^{N-1})}^2 ,
\end{equation}
where $\mathbf{n}_0=\bx/|\bx|$.
\end{proposition}

To prove coercivity, i.e. equation \eqref{eqn:CoercivityEstimSHarpIndef}, it is a 
standard technique (see e.g. 
\cite{Dambrine2011:shapeSensitivity,Dambrine2002:ShapeHessian,Mazari2020:QuantitaiveShrodinger}) 
to expand 
\begin{equation}\label{eq:Xvsvphi}
\mathbf{X}_{\mathcal{A}} (r_0,\theta)\cdot \mathbf{n}_0 = \frac{( r_0+\varphi_{\mathcal{A}}(\theta) )^N - r_0^N}{N \, \rho^{N-1}} 
= \frac{r_0}{N}\left(\left(1+
\frac{\varphi_{\mathcal{A}}(\theta)}{r_0} \right)^N - 1\right) 
\end{equation}
in spherical harmonics: let 
us denote with $(S_k)_k$ the orthonormal basis of $L^2(\sphere^{N-1})$ constituted by the 
eigenfunctions of the Laplace-Beltrami operator $-\Delta_{\sphere^{N-1}}$, with 
eigenvalues $(\sigma_k)_k$:
\begin{equation}\label{eq:spher_arm}
-\Delta_{\sphere^{N-1}} S_k = \sigma_k S_k.
\end{equation}
Then it is well known that $\sigma_0=0$, with constant eigenfunction, while
\[
\sigma_1=\sigma_2=\dots=\sigma_N=N-1,
\] 
and the associated eigenspace is spanned by the orthogonal family of homogeneous harmonic 
polynomials of degree 1 composed by the $N$ cartesian coordinates $x_1, \dots, x_N$; finally, 
$\sigma_k\ge \sigma_{N+1}>N-1$ for $k\ge N+1$. 

A first crucial comment is that, since $\Lcal(\Acal)=1$ and $\bari(\Acal)=0$, in the expansion of 
$\mathbf{X}_{\mathcal{A}} (r_0,\cdot)\cdot \mathbf{n}_0$ the terms corresponding to $S_k$, 
$k\le N$, are negligible.
\begin{lemma}\label{lem:firstspherterms}
Under the assumptions of Proposition \ref{prop:CoercRemainderSharpIndef}, let us write
\begin{equation}\label{eq:exp_spher_X}
\mathbf{X}_{\mathcal{A}}(r_0,\theta) \cdot \mathbf{n}_0 = \sum_{k=0}^{+\infty} c_k S_k(\theta),\qquad 
\text{where }c_k \coloneqq (\mathbf{X}_{\mathcal{A}}(r_0,\cdot) \cdot \mathbf{n}_0, S_k)_{L^2(\sphere^{N-1})}.
\end{equation}
Then $c_0=0$ and 
\begin{equation}\label{eqn:eqn0CoercivitykEstimSHarpIndef}
\sum_{k=N+1}^{+\infty} c_k^2 = \| \mathbf{X}_{\mathcal{A}} \cdot \mathbf{n}_0 \|_{L^2(\partial B)}^2 \left( 1 + \eta \left( \| \varphi_{\mathcal{A}} \|_{L^{\infty}(\sphere^{N-1})} \right) \right) ,
\end{equation}
for some modulus of continuity $\eta$.
\end{lemma}
\begin{proof}
First of all, $c_0 = 0$ by Lemma \ref{lemma:DivFreeSpherDeformSharpIndef} (iii). Next, 
recalling \eqref{eq:Xvsvphi}, let us consider the coefficients
\begin{equation}\label{eqn:eqn1CoercivitykEstimSHarpIndef}
c_k = \int_{\sphere^{N-1}} \theta_k \, \mathbf{X}_{\mathcal{A}} (r_0,\theta)\cdot \mathbf{n}_0 = \frac{r_0}{N} \int_{\sphere^{N-1}} \theta_k \, \left[ ( 1+\varphi_{\mathcal{A}}(\theta)/r_0 )^N - 1 \right] , \quad k=1, \dots, N .
\end{equation}
Since the barycenter of 
$\mathcal{A}$ coincides with $\mathbf{0}$ by assumption, we can write
\begin{equation}\label{eqn:eqn2CoercivitykEstimSHarpIndef}
0 = \int_{\mathcal{A}} x_k - \int_{B} x_k = \frac{r_0^{N+1}}{N+1} \int_{\sphere^{N-1}} \theta_k \,  \left[ ( 1+\varphi_{\mathcal{A}}(\theta)/r_0 )^{N+1} - 1 \right].
\end{equation}
Now, a direct computation gives
\begin{align}\label{eqn:eqn3CoercivitykEstimSHarpIndef}
\begin{split}
&\left \vert \frac{r_0}{N} \int_{\sphere^{N-1}} \theta_k \, \left[ ( 1+\varphi_{\mathcal{A}}(\theta)
/r_0 )^N - 1 \right] - \frac{r_0}{N+1} \int_{\sphere^{N-1}} \theta_k \,  \left[ ( 1+\varphi_{\mathcal{A}}(\theta)/r_0 )^{N+1} - 1 \right] \right \vert \le \\
& \le \eta \left( \| \varphi_{\mathcal{A}} \|_{L^{\infty}(\sphere^{N-1})} \right) \| \varphi_{\mathcal{A}} \|_{L^{2}(\sphere^{N-1})} \le \eta \left( \| \varphi_{\mathcal{A}} \|_{L^{\infty}(\sphere^{N-1})} \right) \| \mathbf{X}_{\mathcal{A}} \cdot \mathbf{n}_0\|_{L^2(
\sphere^{N-1})} ,
\end{split}
\end{align}
where in the last passage we have used \eqref{eqn:EstXvarphiDiffL2SharpIndef}. Hence, combining \eqref{eqn:eqn1CoercivitykEstimSHarpIndef}, \eqref{eqn:eqn2CoercivitykEstimSHarpIndef} and \eqref{eqn:eqn3CoercivitykEstimSHarpIndef} we can deduce \eqref{eqn:eqn0CoercivitykEstimSHarpIndef}.
\end{proof}

To prove  Proposition \ref{prop:CoercRemainderSharpIndef} we will plug the expansion 
\eqref{eq:exp_spher_X} into \eqref{eqn:2ndOrder0EigenExpSharpIndef}. To clarify the role of 
$\dot w$ in the latter, it is convenient to expand in spherical harmonics also the associated 
differential problem \eqref{eqn:eqnDotut0SharpInder}.
\begin{lemma}\label{lem:nomoreBessel}
For every integer $k\ge1$ the problem
\begin{equation}\label{eqn:transmis_for_Sk}
\begin{cases}
-\Delta z_k = \tilde{\lambda}_0 \tilde{m}_0 z_k  & \text{in } \R^N \setminus\partial B, \\
[z_k] = 0, \quad [\partial_{\rho} z_k](\bx) =  j S_k(\bx/|\bx|)  & \text{on } \partial B \\
\int_{\R^N} \tilde m_0 z w =  0 ,
\end{cases}
\end{equation}
where 
\[
[\partial_{\rho} z_k](r_0\theta)=\partial_{\rho}  z_k(r_0^-\theta)-\partial_{\rho} 
z_k(r_0^+\theta) 
\qquad\text{and}\qquad 
j = \tilde\lambda_0(\overline{m}+\underline{m})w(r_0)>0,
\]
has a unique solution $z_k\in H^1(\R^N)$. Moreover:
\begin{enumerate}
\item\label{i:noBess1} $z_k(r\theta)= g_k(r) S_k(\theta)$, where
\[
\begin{cases}
(r^{N-1}g'_k)' + m(r) g_k = r^{N-3} \sigma_k g_k & \text{for } r>0,\ r\neq r_0, \\
g_k(0) = g_k(+\infty) = g_k(r_0^-) - g_k(r_0^+) = 0, \\
g_k'(r_0^-) - g_k'(r_0^+) = j,
\end{cases}
\]
$\sigma_k$ is defined in \eqref{eq:spher_arm}, and $m(r) = \tilde{\lambda}_0 \tilde{m}_0(r\mathbf{e})$, for every $|\mathbf{e}|=1$;
\item\label{i:noBess2} if $k=1,\dots,N$ then $z_k = \partial_{x_k} w$ and $g_k(r) = - w'(r)$;
\item\label{i:noBess3} for every $k\ge1$, $g_k\ge0$ on $(0,+\infty)$;
\item\label{i:noBess4} $g_k(r_0)$ is strictly monotone decreasing with respect to $\sigma_k$; in particular, 
\begin{equation*}
g_k(r_0) \le g_{N+1}(r_0) < g_N(r_0) = -w'(r_0)=\left.-\partial_\rho w\right|_{\partial B},\qquad
\text{ for every $k\ge N+1$.}
\end{equation*}
\end{enumerate}
\end{lemma}
Actually, $g_k$ can be written explicitly in terms of Bessel's functions, yielding their further 
properties. Nonetheless, for the reader's convenience we provide a self-contained proof.  
\begin{proof}
The well-posedness of the problem for $z_k$ follows by the Lax-Milgram theorem, arguing as in  
Lemma \ref{lem:LM}: indeed, notice that the jump condition on $\partial B$ can be 
inserted in the weak formulation as an $H^{-1}$ term $h$ in the right hand side, as it was obtained 
in Lemma \ref{lemma:eqnDotutSharpInder}; in particular, the condition $\langle h, \psi\rangle = 0$ for every $\psi\in H^1(\R^N)$ requires the restriction $k\neq0$.

Once the solution is unique, property \ref{i:noBess1} easily follows by separation of variables, 
using \eqref{eq:spher_arm} and the well-posedness of the problem for $g_k$, and it can be checked 
by direct computation. Uniqueness also implies property \ref{i:noBess2}, since 
\[
\partial_{x_k} w = w'(|\bx|) \cdot\frac{x_k}{|\bx|} = w'(r) S_k(\theta),
\qquad k=1,\dots,N
\]
(with the abuse of notation $w(\bx)=w(|\bx|)=w(r)$) and, in particular, using the equation for $w$,
\[
[w''] = [-\tilde{\lambda}_0 \tilde{m}_0w] = -\tilde\lambda_0(\overline{m}+\underline{m})w(r_0) \qquad
\text{on } \partial B,
\]
so that $g_k=-w'(r)$ and $z_k = \partial_{x_k} w$, $1\le k\le N$. Since $w'(r)<0$ for $r>0$, we also infer property \ref{i:noBess3} for $k\le N$.

In order to prove properties \ref{i:noBess3} and \ref{i:noBess4}, we first need a 
Sturm-Picone-type comparison result for 
solutions with derivative jumps: let $0\le r_1 < r_2 \le + \infty$, $h,k\ge1$, and let us test 
the equation for $g_k$ with $r^{N-1}g_h$ on $(r_1,r_2)$. Then, in case 
$r_0\not\in(r_1,r_2)$, we have
\[
\int_{r_1}^{r_2} r^{N-1}\left(g_h'g_k' - m(r)g_hg_k\right)dr =  
\left[r^{N-1} g_h g'_k\right]_{r_1}^{r_2} - \sigma_k\int_{r_1}^{r_2} r^{N-3}g_hg_k\,dr,
\]
while if $r_0\in(r_1,r_2)$ then we must integrate separately on the two subinterval and 
add the term $[r^{N-1}g_h g'_k]_{r_0^+}^{r_0^-}$ to the right hand side. Exchanging the role of $h$ and $k$, and assuming $g_k (r_1) = g_k (r_2) = 0$, we finally obtain 
\begin{equation}\label{eq:sturm}
(\sigma_k-\sigma_h)\int_{r_1}^{r_2} r^{N-3}g_hg_k\,dr + j(g_k(r_0) - g_h(r_0))\ind{(r_1,r_2)}(r_0) = \left[r^{N-1} 
g_h g'_k\right]_{r_1}^{r_2}.
\end{equation}

Now, to show property \ref{i:noBess1} for $k>N$ let us assume by contradiction that 
$g_k$, $k>N$, is negative somewhere, and denote with $(r_1,r_2)$ a maximal interval of negativity: 
$g_k (r_1) = g_k (r_2) = 0$, $g_k(r)<0$ for $r_1<r<r_2$. Then $g'_k (r_1) \le 0 \le g'_k (r_2)$ 
and, recalling that $g_N=-w'>0$ on $\R_+$, we infer
\[
(\sigma_k-\sigma_N)\int_{r_1}^{r_2} r^{N-3}g_Ng_k\,dr<0,\quad 
j(g_k(r_0) - g_N(r_0))\ind{(r_1,r_2)}(r_0) \le 0,\quad
\left[r^{N-1} g_N g'_k\right]_{r_1}^{r_2}\ge 0,
\]
in contradiction with \eqref{eq:sturm} (with $h=N$), thus property \ref{i:noBess3} follows. 

Finally, take 
any $h<k$ such that $\sigma_h<\sigma_k$. Thanks to property \ref{i:noBess3} we obtain
\[
(\sigma_k-\sigma_h)\int_{0}^{+\infty} r^{N-3}g_hg_k\,dr>0,\qquad 
\left[r^{N-1} g_h g'_k\right]_{0}^{+\infty}= 0,
\]
and \eqref{eq:sturm} (with $r_1=0$, $r_2=+\infty$) implies
\[
j(g_k(r_0) - g_h(r_0))<0,
\]
yielding property \ref{i:noBess4}.
\end{proof}

\begin{proof}[Proof of Proposition \ref{prop:CoercRemainderSharpIndef}]
As a consequence of \eqref{eq:exp_spher_X}, \eqref{eqn:eqn0CoercivitykEstimSHarpIndef}, the proof 
of \eqref{eqn:CoercivityEstimSHarpIndef} can be reduced to the following statement: there exists 
$C=C(\overline{m},\underline{m},N)>0$ such that, for any nearly spherical set $\mathcal{A}$ of 
class $C^{1, 1}$ with $\| \varphi_{\mathcal{A}} \|_{C^{1, 1}(\sphere^{N-1})} < \delta$, we have
\begin{equation}\label{eqn:CoercivitykEstimSHarpIndef}
\left.\ddot \lambda_t\right|_{t=0}\ge C \sum_{k=N+1}^{+\infty} c_k^2 .
\end{equation}

To prove \eqref{eqn:CoercivitykEstimSHarpIndef} we insert  \eqref{eq:exp_spher_X} into 
\eqref{eqn:2ndOrder0EigenExpSharpIndef}. On the one hand, also recalling Lemma 
\ref{lem:nomoreBessel}, \ref{i:noBess2},
\begin{equation}\label{eq:fincoerc1}
\left\vert \frac{\partial w}{\partial \rho} \right\vert_{\partial B}\int_{\partial B} \vert \mathbf{X}_{\mathcal{A}} \vert^2 =
g_N(r_0)\int_{\partial B} \vert \mathbf{X}_{\mathcal{A}} \cdot \mathbf{n}_0 \vert^2 =
g_N(r_0)\sum_{k=1}^{+\infty} c_k^2.
\end{equation}
On the other hand, Lemmas \ref{lemma:eqnDotutSharpInder} and \ref{lem:nomoreBessel}, together 
with the expansion \eqref{eq:exp_spher_X}, yield 
\[
\dot w(r\theta) = \sum_{k=1}^{+\infty} c_k g_k(r) S_k(\theta),
\]
whence
\begin{equation}\label{eq:fincoerc2}
\int_{\partial B} \dot w \mathbf{X}_{\mathcal{A}} \cdot \mathbf{n}_0 
= \sum_{k=1}^{+\infty} g_k(r_0) c_k^2.
\end{equation}

Plugging \eqref{eq:fincoerc1}, \eqref{eq:fincoerc1} into \eqref{eqn:2ndOrder0EigenExpSharpIndef}, 
and using Lemma \ref{lem:nomoreBessel}, \ref{i:noBess2},  
we finally obtain 
\[
\begin{split}
\ddot{\lambda}_t \vert_{t = 0} &= 
2 \tilde \lambda_0 (\overline{m}+\underline{m}) w(r_0) 
\sum_{k=1}^{+\infty} \left(g_N(r_0) - g_k(r_0)\right) c_k^2\\
&=  2 \tilde \lambda_0 (\overline{m}+\underline{m}) w(r_0) 
\sum_{k=N+1}^{+\infty} \left(g_N(r_0) - g_k(r_0)\right) c_k^2,
\end{split}
\]
and \eqref{eqn:CoercivitykEstimSHarpIndef} follows with 
\[
C=2 \tilde \lambda_0 (\overline{m}+\underline{m}) w(r_0)\left[ g_N(r_0) - g_{N+1}(r_0)\right],
\] 
which is strictly positive (and universal) by Lemma \ref{lem:nomoreBessel}, \ref{i:noBess4}.
\end{proof}

\subsection{Continuity of the remainder}\label{sec:cont_of_remaind}

The aim of this section is to prove the following.
\begin{proposition}\label{prop:RemainderrEstimSHarpIndef}
There exist $\delta>0$ sufficiently small and a modulus of continuity $\eta$ such that, for any nearly spherical set $\mathcal{A}$ of class $C^{1, 1}$ with $\| \varphi_{\mathcal{A}} \|_{C^{1, 1}(\sphere^{N-1})} < \delta$ it holds that
\begin{equation}\label{eqn:RemainderrEstimSHarpIndef}
\vert \ddot \lambda_t - \ddot \lambda_0 \vert \le \eta \left( \| \varphi_{\mathcal{A}} \|_{C^{1, 1}(\sphere^{N-1})} \right) \| \mathbf{X}_{\mathcal{A}}(r_0,\cdot) \cdot \mathbf{n} 
\|_{L^2(\sphere^{N-1})}^2
\end{equation}
uniformly for $t \in [0, 1]$.
\end{proposition}
\begin{remark}
Notice that, as it will be clear from the following proofs and as remarked in \cite{Mazari2020:QuantitaiveShrodinger}, Proposition \ref{prop:RemainderrEstimSHarpIndef} holds even for nearly spherical sets $\mathcal{A}$ of class $C^{1, \alpha}$, under the only constraint that $\| \varphi_{\mathcal{A}} \|_{C^{1, \alpha}(\sphere^{N-1})} $ is small enough, for some $0 < \alpha < 1$. Of course, in such a case the $C^{1, 1}$ norm in \eqref{eqn:RemainderrEstimSHarpIndef} has to be substituted by the $C^{1, \alpha}$ norm.
\end{remark}

Before giving the proof, it is useful to introduce two lemmas concerning the $C^{1, \alpha}$ convergence of $u_t$ to $w$ and the $C^{0, \alpha}$ convergence of $\dot u_t$ to $\dot w$, respectively. Since we obtain such results by contradiction, exploiting the convergence induced by elliptic regularity, it is convenient to work in the case $\alpha<1$.  

For any nearly spherical set $\mathcal{A}$ of class $C^{1, \alpha}$, with $0\le \alpha < 1$, we denote with $u_{\mathcal{A}}$ the eigenfunction over the whole $\R^N$ associated to the eigenvalue $\lambda^1(\mathcal{A}, \R^N)$. As already remarked, such function exists and is uniquely defined, up to a normalization that we choose for instance to be in $L^2(\R^N)$. Notice that one could equivalently choose $u_{\mathcal{A}}$ normalized such that $\int_{\R^N} m_{\mathcal{A}} u_{\mathcal{A}}^2 = 1$, and the following discussion would remain unchanged. Then, we proceed with the following

\begin{lemma}\label{lemma:u_tC1aConcSharpIndef}
Fix $0\le \alpha < 1$, $R>0$. There exists a modulus of continuity $\eta$ such that, for any nearly spherical set $\mathcal{A}$ of class $C^{1, \alpha}$, it holds that 
\[
\| u_{\mathcal{A}} - w  \|_{C^{1, \alpha}\left( \overline{B_R} \right)} \le \eta \left( \| \varphi_{\mathcal{A}} \|_{C^{1, \alpha}(\sphere^{N-1})} \right) .
\]
\end{lemma}
\begin{proof}
We proceed by contradiction. Suppose the existence of a sequence of nearly spherical set $\mathcal{A}_k$ of class $C^{1, \alpha}(\partial B)$, satisfying $\| \varphi_{\mathcal{A}} \|_{C^{1, \alpha}(\sphere^{N-1})} \to 0$ for $k \to +\infty$, but such that $\| u_{\mathcal{A}_k} - w  \|_{C^{1, \alpha}\left( \overline{B_R} \right)} \ge \delta$, for any $k \in \N$ and for some $\delta>0$.

However, proceeding as in Section \ref{section:blowUpIndef}, and in particular in Lemma \ref{lem:H1strong}, one gets convergence in $C^{1, \alpha}\left( \overline{B_R} \right)$ of $u_{\mathcal{A}_k}$ to $ w $, both with the same normalization. Hence we have reached a contradiction and the proof is concluded.
\end{proof}

Now we turn to $\dot u_t$. For any nearly spherical set $\mathcal{A}$ consider its flow $\mathcal{A}_t \coloneqq \Phi_{\mathcal{A}}(t, B)$. Denoting with $u_{\mathcal{A}_t}$ the eigenfunction associated to $\mathcal{A}_t$, we can prove the following
\begin{lemma}\label{lemma:dotu_tC0aConcSharpIndef}
Fix $0\le \alpha < 1$, $R>0$. There exists a modulus of continuity $\eta$ such that, for any nearly spherical set $\mathcal{A}$ of class $C^{1, \alpha}(\partial B)$, it holds that 
\[
\| \dot u_{\mathcal{A}_t} - \dot  w   \|_{C^{0, \alpha}\left( \overline{B_R} \right)} \le \eta \left( \| \varphi_{\mathcal{A}_t} \|_{C^{1, \alpha}(\sphere^{N-1})} \right) ,
\]
where $\eta$ is independent of $t \in [0, 1]$.
\end{lemma}
\begin{proof}
Without loss of generality we can assume $B\subset\subset B_R$. First, we notice that $\dot u_{\mathcal{A}_t}$ actually belongs to $C^{0, \alpha}\left( \overline{B_R} \right)$ for any fixed $R>0$ and any $t \in [0, 1]$. To this aim, we consider the function $z_{\mathcal{A}_t} \coloneqq \dot u_{\mathcal{A}_t} + \nabla u_{\mathcal{A}_t} \cdot \mathbf{X}_{\mathcal{A}}$.

Since the function $\nabla u_{\mathcal{A}_t}$ solves a problem analogous to \eqref{eqn:eqn3C2aPropSharpIndef}, using \eqref{lemma:eqnDotutSharpInder} one can prove that the function $z_{\mathcal{A}_t}$ solves
\begin{equation}\label{eqn:eqnDotztSharpInder}
\begin{cases}
-\Delta z_{\mathcal{A}_t} = \dot \lambda_t m_t u_t + \lambda_t m_t z_{\mathcal{A}_t} - \diverg(\nabla (\nabla u_{\mathcal{A}_t} \cdot \mathbf{X}_{\mathcal{A}})) - \lambda_t m_t \nabla u_{\mathcal{A}_t} \cdot \mathbf{X}_{\mathcal{A}}  & \text{in } \R^N , \\
[z_{\mathcal{A}_t}] = 0, \quad [\partial_{\mathbf{n}} z_{\mathcal{A}_t}] =  0  & \text{on } \partial \mathcal{A}_t
\end{cases}
\end{equation}
Since $\nabla u_{\mathcal{A}_t} \in H^{1, p}(B_R)$ for any $1\le p < +\infty$ (this follows by elliptic regularity in $L^p$ spaces for $u_{\mathcal{A}_t}$), standard elliptic regularity for \eqref{eqn:eqnDotztSharpInder} tells us that $z_{\mathcal{A}_t} \in C^{0, \alpha}\left( \overline{B_R} \right)$ for any $\alpha \in [0, 1)$. Moreover, by bootstrap, the bound is uniform with respect to $t \in [0, 1]$ and $\mathcal{A}$, since $\| u_{\mathcal{A}_t} \|_{H^1(B_R)}$ is bounded uniformly with respect to $t \in [0, 1]$ and $\mathcal{A}$,  when $\| \varphi_{\mathcal{A}_t} \|_{C^{1, \alpha}(\sphere^{N-1})}$ is sufficiently small. The last assertion follows from Lemma \ref{lemma:u_tC1aConcSharpIndef}. Consequently, the same regularity property holds for $\dot u_{\mathcal{A}_t}$.

Now we proceed by contradiction. Suppose the existence of a sequence $\mathcal{A}_{t_k}$ of nearly spherical sets of class $C^{1, \alpha}$ with $\| \varphi_{\mathcal{A}_t} \|_{C^{1, \alpha}(\sphere^{N-1})} \to 0$ for $k \to +\infty$, but $\| \dot u_{\mathcal{A}_t} - \dot  w   \|_{C^{0, \alpha}\left( \overline{B_R} \right)} \ge \delta$ for some constant $\delta>0$. We remark that the notation $\mathcal{A}_{t_k}$ is to be intended in the sense that both $\mathcal{A}$ and $t$ can change with $k$. To this sequence we associate the functions $\dot u_{t_k}$ and $\dot w_k$.

Notice that, since for any $0 < \alpha < 1$ the quantity $\|\mathbf{X}_{\mathcal{A}_k}\|_{C^{1, \alpha}(B_{R})}$ is bounded uniformly in $k$, up to a subsequence there exists a limit $\mathbf{X} \in C^{1, \alpha}(B_{R})$.

The first thing we do is proving an $H^1(\R^N)$ bound for $\dot u_{t_k}$ uniformly in $k$. This can be done supposing that $\| \dot u_{t_k} \|_{H^1(\R^N)}$ diverges, and then reaching a contradiction. Since this is a method that we have already used many times, we will be short.

Rescaling \eqref{eqn:eqnDotutTSharpInder} by $\| \dot u_{t_k} \|_{H^1(\R^N )}$, up to a subsequence the functions $\dot u_{t_k}$ converge weakly in $H^1(\R^N)$ and in $C^{0, \alpha}\left( \overline{B_R} \right)$ (thanks to the uniform bound found above) to a solution $s \in H^1(\R^N)$ of
\[
\begin{cases}
-\Delta s = \tilde \lambda_0 \tilde m_0 s  & \text{in } \R^N , \\
[s] = 0, \quad [\partial_{\mathbf{n}} s] =  0  & \text{on } \partial \mathcal{A}_t \\
2 \int_{\R^N} \tilde m_0 s w = 0 .
\end{cases}
\]
Hence $s$ is a multiple of $w$, but the integral condition implies that $s \equiv 0$. Hence, $\dot u_{t_k} \to 0$ in $C^{0, \alpha}\left( \overline{B_R} \right)$. However, from the rescaled version of \eqref{eqn:eqnDotutTSharpInder}, one can notice that for $k$ sufficiently large it holds that $\| \dot u_{t_k} \|_{H^1(\R^n)} \le C \| \dot u_{t_k} \|_{L^2(B_R)} + o(1)$ for some universal constant $C>0$, which leads to a contradiction.

We have just proved that $\| \dot u_{t_k} \|_{H^1(\R^n)}$ is uniformly bounded with respect to $k$. Of course, analogous properties also hold for $\dot w_k$.

To conclude, since $\dot \lambda_{t_k} \to 0$ for $k \to +\infty$ thanks to Lemma \ref{lemma:u_tC1aConcSharpIndef}, and since $u_{t_k} \to w$ strongly in $H^1(\R^N)$ (proceeding as in Sections \ref{section:blowUpIndef} and \ref{section:connectedness}), passing to the limit in \eqref{eqn:eqnDotutTSharpInder} and \eqref{eqn:eqnDotut0SharpInder} it can be noticed that both the sequences $\dot u_{t_k}$ and $\dot w_k$ converge up to a subsequence, weakly in $H^1(\R^N$ and in $C^{0, \alpha}\left( \overline{B_R} \right)$ to a solution in $H^1(\R^N)$ of the problem 
\begin{equation}\label{eqn:eqn1dotu_tC0aConcSharpIndef}
\begin{cases}
-\Delta v = \tilde \lambda_0 \tilde m_0 v  & \text{in } \R^N , \\
[v] = 0, \quad [\partial_{\mathbf{n}} v] =  \tilde \lambda_0 (\overline{m}+\underline{m}) w \, \mathbf{X} \cdot \mathbf{n}  & \text{on } \partial B \\
2 \int_{\R^N} \tilde m_0 v w = - (\overline{m} + \underline{m}) \int_{\partial B} w^2 \, \mathbf{X} \cdot \mathbf{n} .
\end{cases}
\end{equation}
However, as already remarked, there exist at most one solution in $H^1(\R^N)$ to problem \eqref{eqn:eqn1dotu_tC0aConcSharpIndef}, so that $\dot u_{t_k}$ and $\dot w_k$ converge in $C^{0, \alpha}\left( \overline{B_R} \right)$ to the same function, and we have reached a contradiction. 
\end{proof}

Now we are ready to prove Proposition \ref{prop:RemainderrEstimSHarpIndef}.
\begin{proof}[Proof of Proposition \ref{prop:RemainderrEstimSHarpIndef}.]
We write $\ddot \lambda_t - \ddot \lambda_0 \coloneqq A - B - C$, with
\begin{align*}
    & A \coloneqq \lambda_t \left[ (\overline{m}+\underline{m}) \int_{\mathcal{A}_t} \diverg(u_t^2 \mathbf{X}_{\mathcal{A}}) \right]^2, \\
    & B \coloneqq 2 (\overline{m}+\underline{m}) \left[ \lambda_t \int_{\partial \mathcal{A}_t} u_t \nabla u_t \cdot \, \mathbf{X}_{\mathcal{A}} \, \mathbf{X}_{\mathcal{A}} \cdot \mathbf{n}_t - \lambda_0 \int_{\partial B} w \nabla w \cdot \mathbf{n} \vert \mathbf{X}_{\mathcal{A}} \vert^2 \right] , \\
    & C \coloneqq 2 (\overline{m}+\underline{m}) \left[ \lambda_t \int_{\partial \mathcal{A}_t} \dot u_t u_t \mathbf{X}_{\mathcal{A}} \cdot \mathbf{n}_t - \lambda_0 \int_{\partial B} \dot w w \mathbf{X}_{\mathcal{A}} \cdot \mathbf{n} \right] .
\end{align*}
Rewriting $A$ as a surface integral on $\partial \mathcal{A}_t$, transporting all the integrals on $\partial B$ and denoting with $J_{T, t}$ the tangential jacobian associated with $\partial \mathcal{A}_t$, the above terms can be rewritten as
\begin{align*}
    & A \coloneqq \lambda_t \left[ (\overline{m}+\underline{m}) \int_{\partial B} (u_t^2 \mathbf{X}_{\mathcal{A}} \cdot \mathbf{n}_t) \circ \Phi_{\mathcal{A}} J_{T, t} \right]^2, \\
    & B \coloneqq 2 (\overline{m}+\underline{m}) \left[ \lambda_t \int_{\partial B} \left( u_t \nabla u_t \cdot \, \mathbf{X}_{\mathcal{A}} \, \mathbf{X}_{\mathcal{A}} \cdot \mathbf{n}_t \right) \circ \Phi_{\mathcal{A}} J_{T, t} - \lambda_0 \int_{\partial B} w \nabla w \cdot \mathbf{n} \vert \mathbf{X}_{\mathcal{A}} \vert^2 \right] , \\
    & C \coloneqq 2 (\overline{m}+\underline{m}) \left[ \lambda_t \int_{\partial B} \left( \dot u_t u_t \mathbf{X}_{\mathcal{A}} \cdot \mathbf{n}_t \right) \circ \Phi_{\mathcal{A}} J_{T, t} - \lambda_0 \int_{\partial B} \dot w w \mathbf{X}_{\mathcal{A}} \cdot \mathbf{n} \right] .
\end{align*}
Taking into account that the $C^{1,\alpha}$ norm of $\varphi_{\Acal}$ is controlled by the $C^{1,1}$ norm, the proof can be concluded using Lemma \ref{lemma:u_tC1aConcSharpIndef} to treat $u_t$, Lemma \ref{lemma:dotu_tC0aConcSharpIndef} to treat $\dot u_t$, Lemma \ref{lemma:DivFreeSpherDeformSharpIndef} (iv), (v) to deal with the composition by $\Phi_{\mathcal{A}}$, Lemma \ref{lemma:GeometricC1aConvSharpIndef} (iii), (iv) (notice that Lemma \ref{lemma:GeometricC1aConvSharpIndef} can be easily rewritten in terms of a modulus of continuity) to deal with $\mathbf{n}_t$ and $J_{T, t}$ and, to conclude, the continuity of $\lambda_t$ stated in Proposition \ref{prop:W2inftyRegulSharpPos}.
\end{proof}

\begin{proof}[End of the proof of Theorem \ref{thm:quantitStabRNSharpPos}.] Expanding $\lambda_t$ in Taylor series up to the second order, we obtain
\[
\lambda^1(\mathcal{A}, \R^N) - \lambda^1(B, \R^N) = \ddot \lambda_0 + \int_{0}^1 (1-t) (\ddot \lambda_t - \ddot \lambda_0) + \int_{0}^1 (1-t) \ddot \lambda_0.
\]
The proof can be concluded using Propositions \ref{prop:RemainderrEstimSHarpIndef}, \ref{prop:CoercRemainderSharpIndef} and Lemma \ref{lemma:DivFreeSpherDeformSharpIndef} (vi) above.
\end{proof}

\textbf{Data Availability Statement.} Data sharing is not applicable to this article as no new data were created or 
analyzed in this study.

\textbf{Acknowledgements.} LF is partially supported by the European Research Council (ERC), under the European Union's Horizon 2020 research and innovation program, through the project ERC VAREG - {\em Variational approach to the regularity of the free boundaries} (grant agreement No. 853404). GV is partially supported by the Portuguese 
government through FCT/Portugal under the project PTDC/MAT-PUR/1788/2020. The authors are members 
of the INdAM-GNAMPA group (``Gruppo Nazionale per l'Analisi Matematica, la Probabilit\`a e le loro 
Applicazioni -- Istituto Nazionale di Alta Matematica'').


\bibliographystyle{abbrv}
\bibliography{Paper2_bib.bib}

\medskip
\small
\begin{flushright}
\noindent 
\verb"lorenzo.ferreri@sns.it"\\
Classe di Scienze, Scuola Normale Superiore\\
piazza dei Cavalieri 7, 56126 Pisa (Italy)\\
\bigskip

\verb"gianmaria.verzini@polimi.it"\\
Dipartimento di Matematica, Politecnico di Milano\\ 
piazza Leonardo da Vinci 32, 20133 Milano (Italy)
\end{flushright}

\end{document}